\documentclass{article}
\font\tenmath=msbm10
\font\sevenmath=msbm7
\font\fivemath=msbm5
\newfam\mathfam \textfont\mathfam=\tenmath
\scriptfont\mathfam=\sevenmath \scriptscriptfont\mathfam=\fivemath
\usepackage{amsmath}

\usepackage{amsfonts}
\usepackage{amssymb}
\usepackage[american]{babel}
\usepackage{graphicx}
\usepackage{flafter}
\usepackage[section]{placeins}
\usepackage{psfrag}
\usepackage{subfigure}

\newcommand{\DS}{\displaystyle}
\newcommand{\field}[1]{\mathbb{#1}}

\newcommand{\N}{{\rm I}\kern-0.18em{\rm N}}

\newcommand{\R}{{\rm I}\kern-0.18em{\rm R}}
\newcommand{\h}{{\rm I}\kern-0.18em{\rm H}}
\newcommand{\K}{{\rm I}\kern-0.18em{\rm K}}
\newcommand{\p}{{\rm I}\kern-0.18em{\rm P}}
\newcommand{\E}{{\rm I}\kern-0.18em{\rm E}}
\newcommand{\Z}{{\rm Z}\kern-0.18em{\rm Z}}

\newcommand{\1}{{\rm 1}\kern-0.24em{\rm I}}

\newcommand{\X}{\field{X}}
\newcommand{\n}{\mathcal{N}}

\newcommand{\ud}{\mathrm{d}}

\newcommand{\epr}{\hfill\hbox{\hskip 4pt\vrule width 5pt
                  height 6pt depth 1.5pt}\vspace{0.5cm}\par}

\newcommand{\e}{\textrm{e}}
\newcommand{\pp}{{\sf p}}
\newcommand{\ubar}{\underbar}
\newcommand{\argmin}{\mathop{\mathrm{argmin}}}
\newcommand{\cI}{\mathcal{I}}
\newtheorem{TH1}{Theorem}[section]

\newtheorem{cor}{Corollary}[section]

\newtheorem{lem}{Lemma}[section]

\newtheorem{defin}{Definition}[section]

\begin{document}
\title{\bf Linear and convex aggregation of density estimators}
\author{\textsc{
Philippe Rigollet \hspace{2cm}  Alexandre B. Tsybakov}\\
\\
\normalsize Laboratoire de Probabilit\'es et Mod\`eles Al\'eatoires,\\
\normalsize Universit\'e Paris 6,\\
\normalsize 4 pl. Jussieu, 75252 Paris Cedex 05, France,\\
\normalsize \{rigollet, tsybakov\}@ccr.jussieu.fr}
\date{\normalsize \today}
\maketitle

\begin{abstract}
We study the problem of linear and convex aggregation of $M$
estimators of a density with respect to the mean squared risk. We
provide procedures for linear and convex aggregation and we prove
oracle inequalities for their risks. We also obtain lower bounds
showing that these procedures are rate optimal in a minimax sense.
As an example, we apply general results to aggregation of
multivariate kernel density estimators with different bandwidths. We
show that linear and convex aggregates mimic the kernel oracles in
asymptotically exact sense for a large class of kernels including
Gaussian, Silverman's and Pinsker's ones. We prove that, for
Pinsker's kernel, the proposed aggregates are sharp asymptotically
minimax simultaneously over a large scale of Sobolev classes of
densities. Finally, we provide simulations demonstrating performance
of the convex aggregation procedure.
\end{abstract}

\medskip

{\footnotesize



{\it Short title:} Aggregation of density estimators.

}

\section{Introduction}
\setcounter{equation}{0} Consider i.i.d. random vectors $X_1,
\ldots, X_n$ with values in $\R^d$ having an unknown common
probability density  $p \in L_2(\R^d)$ that we want to estimate. For
an estimator $\hat p$ of $p$ based on the sample ${\X}^n=(X_1,
\ldots, X_n)$, define the $L_2$-risk
$$
R_n(\hat p,p) = E_p^n\|{\hat p} - p\|^2,
$$
where $E_p^n$ denotes the expectation w.r.t. the distribution
$P_p^n$ of ${\X}^n$ and, for a function $g \in L_2(\R^d)$,
$$
\|g\| = \left(\int_{\R^d} g^2(x) \ud x \right)^{1/2}.
$$
Suppose that we have $M\ge 2$ estimators $\hat p_{1}, \dots, \hat
p_{M}$ of the density $p$ based on the sample ${\X}^n$. The problem
that we study here is to construct a new estimator $\tilde p_n$ of
$p$, called {\it aggregate}, which is approximately at least as good
as the best linear or convex combination of $\hat p_{1}, \dots, \hat
p_{M}$. The problems of linear and convex aggregation of density
estimators under the $L_2$ loss can be stated as follows.

\begin{enumerate}
\item {\bf Problem (L): linear aggregation.}
Find a \emph{linear aggregate}, i.e. an estimator $\tilde p_n^{\bf
L}$ which satisfies
\begin{equation}
\label{int1}
 R_n(\tilde p_n^{\bf L},p) \le \inf_{\lambda \in
\R^M}R_n({\sf p}_\lambda,p) + \Delta_{n,M}^{\bf L}
\end{equation}
for every $p$ belonging to a large class of densities $\mathcal{P}$,
where
$$
{\sf p}_\lambda =  \sum_{j=1}^M \lambda_j \hat p_{j}, \quad \lambda
= (\lambda_1,\dots,\lambda_M),
$$
and $\Delta_{n,M}^{\bf L}$ is a sufficiently small remainder term
that does not depend on $p$.

\item {\bf Problem (C): convex aggregation.}
Find a \emph{convex aggregate}, i.e. an estimator $\tilde p_n^{\bf
C}$ which satisfies
\begin{equation}
\label{int2}
 R_n(\tilde p_n^{\bf C},p) \le \inf_{\lambda \in H} R_n({\sf p}_\lambda,p) + \Delta_{n,M}^{\bf C}
\end{equation}
for every $p$ belonging to a large class of densities $\mathcal{P}$,
where
 $\Delta_{n,M}^{\bf C}$ is a sufficiently small remainder term that
 does not depend on $p$, and $H$ is a convex compact subset of
 $\R^M$. We will discuss in more detail the case $H=\Lambda^M$
 where $\Lambda^M$ is a simplex,
$$
\Lambda^M = \Big\{ \lambda\in\R^M: \ \lambda_j\ge 0,  \sum_{j=1}^M
\lambda_j
  \le 1\Big\}.
$$
\end{enumerate}
Our aim is to find aggregates satisfying \eqref{int1} or
\eqref{int2} with the smallest possible remainder terms
$\Delta_{n,M}^{\bf L}$ and $\Delta_{n,M}^{\bf C}$. These remainder
terms characterize the price to pay for aggregation.

Linear and convex aggregates mimic the best linear (respectively,
convex) combinations of the initial estimators. Along with them, one
may consider {\it model selection (MS) aggregates} that mimic the
best among the initial estimators $\hat p_{1}, \dots, \hat p_{M}$.
We do not analyze this type of aggregation here.

The study of convergence properties of aggregation methods has been
initiated by Nemirovski (2000), Catoni (1999, 2004) and Yang (2000).
Most of the results were obtained for the regression and Gaussian
white noise models (see a recent overview in Bunea, Tsybakov and
Wegkamp (2004)). Aggregation of density estimators has received less
attention. The work on this subject is mainly devoted to the MS
aggregation with the Kullback-Leibler divergence as a loss function
[Catoni (1999, 2004), Yang (2000),  Zhang (2003)], and is based on
information-theoretical ideas close to the earlier papers of Barron
(1987), Li and Barron (1999). Devroye and Lugosi (2001) developed a
method of MS aggregation of density estimators satisfying certain
complexity assumptions under the $L_1$ loss.

To our knowledge, linear aggregation of density estimators has not
been previously studied. For convex aggregation, the only paper we
are aware of is that of Birg\'e (2003) where this type of
aggregation under the $L_1$ loss is considered, while we study here
the $L_2$ loss. In his setup, Birg\'e (2003) proves an inequality
which is weaker than \eqref{int2}, with the oracle risk on the right
hand side multiplied by a constant which is much larger than 1.

We do not only suggest aggregates satisfying sharp oracle
inequalities \eqref{int1}, \eqref{int2}, but also demonstrate their
optimality. Namely, we introduce the notion of optimal rate of
aggregation and show that our aggregates attain optimal rates. This
extends to density estimation context some results of the paper of
Tsybakov (2003) where optimal rates of aggregation for the
regression model have been obtained.

The main purpose of aggregation is to improve upon the initial set
of estimators $\hat p_{1}, \dots, \hat p_{M}$. This is a general
tool that applies to various kinds of estimators satisfying very
mild conditions (we only assume that they are square integrable).
Consider, for example, the simplest case when we have only two
estimators ($M=2$), where $\hat p_{1}$ is a good parametric density
estimator for some fixed regular parametric family and $\hat p_{2}$
is a nonparametric density estimator. If the underlying density $p$
belongs to the parametric family, $\hat p_{1}$ is perfect: its risk
converges with the parametric rate $O(1/n)$. But for densities which
are not in this family it may not converge at all. As for $\hat
p_{2}$, it converges with a slow nonparametric rate even if the
underlying density is within the parametric family. Aggregation (cf.
Section 2 below) allows one to construct procedures that combine the
advantages of both $\hat p_{1}$ and $\hat p_{2}$: the convex or
linear aggregates converge with the parametric rate $O(1/n)$ if $p$
is within the parametric family, and with a nonparametric rate
otherwise. Similar use of aggregation can be done in the problem of
adaptation to the unknown smoothness (cf. Sections 5 and 6). In this
case the index $j$ of $\hat p_{j}$ corresponds to a value of the
smoothing parameter, and the adaptive estimators in the oracle or
minimax sense can be obtained as linear or convex aggregates. Of
course, there exists a large variety of other methods of adaptation
to unknown smoothness. In the numerical examples that we consider,
our aggregates are comparable to benchmarks, and show somewhat more
stable behavior for densities with highly inhomogeneous smoothness
(cf. Section 7). It is important to note that aggregation can be
used for adaptation to other characteristics than smoothness, for
example, to the dimension of the subspace where the data effectively
lie, under dimension reduction models [cf. Samarov and Tsybakov
(2005)].

In this paper, we consider only one example of application of our
general results to the problem of adaptation to the unknown
smoothness. Specifically, we deal with aggregation of multivariate
kernel density estimators with different bandwidths. Here the number
$M=M_n$ of the estimators depends on $n$ and satisfies $M_n/n \to
0$, as $n \to \infty$. We show in Corollary \ref{corkernoracle} that
linear and convex aggregates mimic the kernel oracles in sharp
asymptotic sense. This corollary is in the spirit of Stone's (1984)
theorem on asymptotic optimality of cross-validation, but it is more
powerful in several aspects because it is obtained under weaker
conditions on $p$ and covers kernels with unbounded support
including Gaussian, Silverman's and Pinsker's kernels. Another
application of our results is that, for Pinsker's kernel, we
construct aggregates that are sharp asymptotically minimax
simultaneously over a large scale of Sobolev classes of densities in
the multidimensional case.

To perform aggregation, we use a sample splitting scheme. The sample
$\X^n$ is split into two independent subsamples $\X^m_1$ (training
sample) and $\X_2^{\ell}$ (validation sample) of sizes $m$ and
${\ell}$ respectively where $m+{\ell}=n$ and usually $m\gg {\ell}$.
The first
 subsample $\X^m_1$ is used to
construct estimators $\hat p_{j} = \hat p_{m,j}, \ j=1, \dots, M$,
while the second subsample $\X^{\ell}_2$ is used to aggregate them,
i.e., to construct ${\tilde p}_n$ (thus, ${\tilde p}_n$ is
measurable w.r.t. the whole sample $\X^n$). In a first analysis we
will not consider sample splitting schemes but rather deal with a
``pure aggregation" framework (as in most of the papers on the
subject, cf. ,e.g., Nemirovski (2000), Juditsky and Nemirovski
(2000) and Tsybakov (2003) for the regression problem) where the
first subsample is frozen. This means that instead of the estimators
$\hat p_{1}, \dots, \hat p_{M}$ we have fixed functions $p_1, \dots,
p_M$ and that the expectations in oracle inequalities are taken
only w.r.t. the second subsample.


This paper is organized as follows. In Section~\ref{sec2} we
introduce linear and convex aggregation procedures and prove that
they satisfy oracle inequalities of the type \eqref{int1} and
\eqref{int2}. Section~\ref{sec3} provides lower bounds showing
optimality of the rates obtained in Section~\ref{sec2}. Consequences
for averaged aggregates are stated in Section~\ref{sec4}. In
Sections~\ref{sec5} and~\ref{sec6} we apply the results of Sections~\ref{sec2} and~\ref{sec4} to aggregation of kernel density
estimators. Section~\ref{sec7} contains a simulation study.
Throughout the paper we denote by $c_i$ finite positive constants.

\section{Oracle inequalities for linear and convex aggregates}
\setcounter{equation}{0} \label{sec2} In this section, $p_1, \ldots,
p_M$ are fixed functions, not necessarily probability densities.
>From now on the notation $\pp_\lambda$ for a vector $\lambda=
(\lambda_1,\dots,\lambda_M)\in \R^M$ is understood in the following
sense:
$$
\pp_\lambda \triangleq \sum_{j=1}^M \lambda_j p_{j},
$$
and, since for any fixed $\lambda \in \R^M$, the function $\pp_\lambda$ is non-random, we have
$$R_n({\sf p}_\lambda,p) = \|\pp_\lambda - p\|^2.
$$
Denote by $\mathcal{P}_0$ the class of all densities on $\R^d$
bounded by a constant $L>0$:
$$
\mathcal{P}_0 \triangleq  \left\{p:\R^d\to \R \, \Big| \, p\ge 0,
\int_{\R^d}p(x) \ud x =1, \, \|p\|_\infty \le
  L\right\},
$$
where $\|\cdot\|_\infty$ stands for the $L_\infty(\R^d)$ norm. The
constant $L$ need not be known to the statistician.

We first give an oracle inequality for linear aggregation. Denote by
$\mathcal{L}$ the linear span of $p_1,\ldots,p_M$. Let
$\phi_1,\ldots,\phi_{M'}$ with $M'\le M$ be an orthonormal basis of
$\mathcal{L}$ in $L_2(\R^d)$. Define a linear aggregate
\begin{equation}
\label{11} {\tilde p}_n^{\bf L}(x) \triangleq \sum_{j=1}^{M'} {\hat
\lambda}_j^{\bf L} \phi_j(x), \quad x\in \R^d,
\end{equation}
where
$$
{\hat \lambda}_j^{\bf L} = \frac{1}{n}\sum_{i=1}^{n} \phi_j(X_i).
$$
\begin{TH1}
\label{t2l} Assume that $p_1,\ldots,p_M\in L_2(\R^d) $ and $p \in
\mathcal{P}_0$. Then
\begin{equation}
\label{bslin} R_n({\tilde p}_n^{\bf L},p)\le \min_{\lambda\in\R^M}
\|\pp_\lambda - p\|^2 + \frac{LM}{{n}}
\end{equation}
for any integers $M\ge 2$ and $n\ge 1$.
\end{TH1}
{\sc Proof.} Consider the projection of $p$ onto $\mathcal{L}$:
$$
p_{\mathcal{L}}^{*} \triangleq  \argmin_{\pp_{\lambda}\in
\mathcal{L}} \|\pp_\lambda - p\|^2 = \sum_{j=1}^{M'}
{\lambda}^*_j\phi_j,
$$
where  ${\lambda}^*_j = (p,\phi_j)$, and $(\cdot,\cdot)$ is the
scalar product in $L_2(\R^d)$. Using the Pythagorean theorem we get
that, almost surely,
$$
\|{\tilde p}_n^{\bf L} -p\|^2 = \sum_{j=1}^{M'} ({\hat
\lambda}_j^{\bf L} -{\lambda}^*_j)^2 + \|p^*_{\mathcal{L}} - p\|^2.
$$
To finish the proof it suffices to take expectations in the last
equation and to note that $E_p^n({\hat \lambda}_j^{\bf L}) =
{\lambda}^*_j$ and
$$
E_p^n\Big[({\hat \lambda}_j^{\bf L} -{\lambda}^*_j)^2\Big] ={\rm
Var}( {\hat \lambda}_j^{\bf L}) \le \frac{1}{n}\int_{\R^d}
\phi_j^2(x) p(x) \ud x \leq \frac{L}{n} \ .
$$
\epr

Consider now convex aggregation. Its aim is to mimic the
\emph{convex oracle} defined as $\lambda^* = \argmin_{\lambda \in
   H}\|\pp_{\lambda}-p\|^2$ where $H$ is a given convex compact subset of
  $\R^M$.
  Clearly,
$$
\|\pp_{\lambda}-p\|^2= \|\pp_{\lambda}\|^2 -2\int_{\R^d}
\pp_{\lambda}p+\|p\|^2.
$$
Removing here the term $\|p\|^2$ independent of $\lambda$ and
estimating $\int_{\R^d} p_j p$ by $n^{-1} \sum_{i=1}^n p_j(X_i)$ we
get the following estimate of the oracle
\begin{equation}
\label{deflamconv} \hat \lambda^{\bf C} = \argmin_{\lambda \in H}
\left\{\|\pp_{\lambda}\|^2 - \frac{2}{n} \sum_{i=1}^n
\pp_{\lambda}(X_i)\right\}.
\end{equation}
Now, we define a \emph{convex aggregate} $\tilde p_n^{\bf C}$ by
$$
\tilde p_n^{\bf C} \triangleq \sum_{j=1}^M\hat \lambda_j^{\bf C}
p_j=\pp_{\hat
  \lambda^{\bf C}} .
$$
\begin{TH1}
\label{oracleconv} Let $H$ be a convex compact subset of $\R^M$.
Assume that $p_1,\ldots,p_M\in L_2(\R^d) $ and $p \in
\mathcal{P}_0$. Then the convex aggregate ${\tilde p}_n^{\bf C}$
satisfies
\begin{equation}
\label{bb} R_n({\tilde p}_n^{\bf C},p)\le \min_{\lambda\in H}
\|\pp_\lambda - p\|^2 + \frac{4LM}{{n}}
\end{equation}
for any integers $M\ge 2$ and $n\ge 1$.
\end{TH1}
{\sc Proof.} We will write for brevity $\hat\lambda=
\hat\lambda^{\bf C}$. First note that the mapping
$\lambda\mapsto\|\pp_{\lambda}\|^2 - \frac{2}{n} \sum_{i=1}^n
\pp_{\lambda}(X_i)$ is continuous, thus $\hat \lambda$ exists, and
the oracle $\lambda^* = \argmin_{\lambda \in
   H}\|\pp_{\lambda}-p\|^2$ also exists.
The definition of ${\hat\lambda}$ implies that, for any $p \in
\mathcal{P}_0$,
\begin{equation}
\label{c1} \|\pp_{\hat\lambda}-p\|^2 \le  \|\pp_{\lambda^*} - p\|^2
+ 2T_n
\end{equation}
where
$$
T_n \triangleq  \frac{1}{n} \sum_{i=1}^n
  \pp_{\hat\lambda-\lambda^*}(X_i)-
  \int_{\R^d} \pp_{\hat\lambda-\lambda^*} p .
$$
Introduce the notation
$$
\mathcal{Z}_{n} \triangleq \sup_{\mu\in \R^M: \, \|\pp_\mu\|\neq 0}
\frac{ \left|\frac{1}{n} \sum_{i=1}^n
  \pp_\mu(X_i)-E_p^n [ \pp_\mu(X_1) ] \right|}
{\|\pp_\mu\|}.
$$
Using the Cauchy-Schwarz inequality, the identity
$\pp_{\hat\lambda-\lambda^*} = \pp_{\hat\lambda}-\pp_{{\lambda}^*}$
and the elementary inequality $2\sqrt{xy}\le ax + y/a,\ \forall \,
x,y,a>0$, we get
\begin{eqnarray}
\nonumber
 E_p^n|T_n|  &\le & E_p^n\left( \mathcal{Z}_{n}
  \|\pp_{\hat\lambda-\lambda^*}\|\right)
  \\
  &\le&  \sqrt{E_p^n( \mathcal{Z}_{n}^2)}
  \sqrt{E_p^n(\|\pp_{\hat\lambda-\lambda^*}\|^2)}\nonumber\\
   &\le& \frac{a}{2}E_p^n(\|\pp_{\hat\lambda}-\pp_{\lambda^*}\|^2)
   +\frac{1}{2a}E_p^n( \mathcal{Z}_{n}^2),\quad \forall \ a>0.\label{c2}
\end{eqnarray}
Representing $\pp_\mu$ in the form
$\pp_\mu=\sum_{l=1}^{M'}\nu_l\phi_l$ where $\nu_l\in \R$ and
$\{\phi_l\}$ is an orthonormal basis in $\mathcal L$ (cf. proof of
Theorem~\ref{t2l}) we find
 $$
   \mathcal{Z}_{n} \le  \sup_{\nu\in\R^M\setminus \{0\}}
\frac{|\sum_{l=1}^{M'}\nu_l\zeta_l|}{|\nu|} =
\Big(\sum_{l=1}^{M'}\zeta_l^2\Big)^{1/2},
$$
where $|\nu|=\Big(\sum_{l=1}^{M'}\nu_l^2\Big)^{1/2}$ and
$$
\zeta_l = \frac{1}{n} \sum_{i=1}^n
  \phi_l(X_i)-E_p^n [ \phi_l(X_1) ].
$$
Hence
\begin{equation}\label{c3}
E_p^n( \mathcal{Z}_{n}^2) \le \frac{M'}{n}\max_{l=1,\dots,M'}E_p^n [
\phi_l^2(X_1) ] \le \frac{LM}{n},
\end{equation}
whenever $\|p\|_\infty \le L$. Since $\{\pp_\lambda: \lambda\in H\}$
is a convex subset of $L_2(\R^d)$ and $\pp_{\lambda^*}$ is the
projection of $p$ onto this set, we have
\begin{equation}\label{c4}
\|\pp_{\lambda}-p\|^2\ge
\|\pp_{\lambda^*}-p\|^2+\|\pp_{\lambda}-\pp_{\lambda^*}\|^2, \quad
\forall \ \lambda\in H, \ p \in L_2(\R^d).
\end{equation}
Using \eqref{c4} with $\lambda=\hat\lambda$, \eqref{c2} and
\eqref{c3} we obtain
$$
E_p^n|T_n| \le
\frac{a}{2}\left\{E_p^n(\|\pp_{\hat\lambda}-p\|^2-\|\pp_{\lambda^*}-p\|^2)
\right\}
   +\frac{LM}{2an} \ .
$$
This and \eqref{c1} yield that, for any $0<a<1$,
$$
E_p^n(\|\pp_{\hat\lambda}-p\|^2) \le \|\pp_{\lambda^*}-p\|^2 +
\frac{LM}{a(1-a)n} \ .
$$
Now, \eqref{bb} follows by taking the infimum of the right hand side
of this inequality over $0<a<1$. \epr

\section{Lower bounds and optimal aggregation}
\setcounter{equation}{0}
\label{sec3}

We first define the notion of \emph{optimal rate of aggregation} for
density estimation, similar to that for the regression problem given
in Tsybakov (2003). It is related to the minimax behavior of the
excess risk
$$
{\mathcal E}({\tilde p}_n,p,H) = R_n({\tilde p}_n,p) -
\inf_{\lambda\in H} \|\pp_\lambda - p\|^2
$$
for a given class $H$ of weights $\lambda$.
\begin{defin}
\label{d1} Let $\mathcal{P}$ be a given class of probability
densities on $\R^d$, and let $H\subseteq \R^M$ be a given class of
weights. A sequence of positive numbers $\psi_{n}(M)$ is called {\bf
optimal rate of
  aggregation} for $H$ over $\mathcal{P}$ if
\begin{itemize}
\item for any functions
$p_j \in L_2(\R^d), j=1,\dots,M,$ there exists an estimator ${\tilde
p}_n$ of $p$ (aggregate) such that
\begin{equation}
\label{defub} \sup_{p\in \mathcal{P}} \left[R_n({\tilde p}_n,p) -
\inf_{\lambda\in H} \|\pp_\lambda - p\|^2\right]
 \le C\psi_{n}(M),
\end{equation}
for any integer $n\ge 1$ and for some constant $C<\infty$
independent of $M$ and $n$,
\end{itemize}
and
\begin{itemize}
\item there exist functions
$p_j \in L_2(\R^d), j=1,\dots,M,$ such that for all estimators $T_n$
of $p$, we have
\begin{equation}
\label{deflb} \sup_{p\in \mathcal{P}} \left[ R_n(T_n,p) -
\inf_{\lambda\in H} \|\pp_\lambda - p\|^2 \right] \ge c\psi_{n}(M),
\end{equation}
for any integer $n\ge 1$ and for some constant $c>0$ independent of
$M$ and $n$.
\end{itemize}
When \eqref{deflb} holds, an aggregate ${\tilde
p}_n$ satisfying \eqref{defub} is called {\bf rate optimal aggregate} for $H$ over
$\mathcal{P}$.
\end{defin}
Note that this definition applies to aggregation of any functions
$p_j$ in $L_2(\R^d)$, they are not necessarily supposed to be
probability densities.

Theorems~\ref{t2l} and~\ref{oracleconv} provide upper bounds of the type
\eqref{defub} with the rate $\psi_n(M) = LM/n$ for linear and convex
aggregates ${\tilde p}_n= {\tilde p}_n^{\bf L}$ and ${\tilde p}_n=
{\tilde p}_n^{\bf C}$ when $\mathcal{P} = \mathcal{P}_0$ and $H=
\R^M$ or $H$ is a convex compact subset of $\R^M$. In this section
we complement these results by lower bounds of the type \eqref{deflb}
showing that $\psi_n(M) = LM/n$ is optimal rate of linear and convex
aggregation. The proofs will be based on the following lemma which
is adapted from Corollary~4.1 of Birg\'e (1986), p. 281.

\begin{lem}
\label{lem1} Let $\mathcal{C}$ be a set of functions of the
following type
$$
\mathcal{C}=\bigg\{f+\sum_{i=1}^r \delta_i g_i, \ \delta_i \in
\{0,1\}, \ i=1, \ldots, r \bigg\},
$$
where the $g_i$ are functions on $\R^d$ with disjoint supports, such
that $\int g_i(x) \ud x =0$, $f$ is a probability density on $\R^d$
which is constant on the union of the supports of $g_i$'s, and
$f+g_i\ge 0$ for all $i$. Assume that
\begin{equation}
\label{condlem1} \min_{1\le i\le r}\| g_i \|^2 \ge \alpha>0 \quad
\textrm{and} \quad \max_{1\le i\le r} h^2(f,f+g_i) \le \beta <1 \,,
\end{equation}
where $h^2(f,g)=(1/2)\int (\sqrt{f(x)}-\sqrt{g(x)})^2 dx$ is the
squared Hellinger distance between two probability densities $f$ and
$g$. Then
$$
\inf_{T_n} \sup_{p \in {\cal C}}R_n(T_n,p) \ge \frac{r \alpha}{4}
(1-\sqrt{2n\beta})
$$
where $\DS \inf_{T_n}$ denotes the infimum over all estimators.
\end{lem}

Consider first a lower bound for linear aggregation of density
estimators. We are going to prove \eqref{deflb} with $\psi_n(M) =
LM/n$, $\mathcal{P} = \mathcal{P}_0$ and $H= \R^M$. Note first that
for $\mathcal{P} = \mathcal{P}_0$ there is a natural limitation on
the value $c\psi_n(M)$ on the right hand side of \eqref{deflb},
whatever is $H$. In fact, $\inf_{T_n} \sup_{p\in \mathcal{P}_0}
\Big[R_n(T_n,p)- \inf_{\lambda\in H} \|\pp_{\lambda} - p\|^2\Big]
\le \inf_{T_n} \sup_{p\in \mathcal{P}_0}R_n(T_n,p)\le \sup_{p\in
\mathcal{P}_0}R_n(0,p) = \sup_{p\in \mathcal{P}_0} \|p\|^2\le L.$
Therefore, we must have $c\psi_n(M)\le L$ where $c$ is the constant
in \eqref{deflb}. For $\psi_n(M) = LM/n$ this means that only the
values $M$ such that $M \le c_0 n$ are allowed, where $c_0>0$ is a
constant. The upper bounds of Theorems~\ref{t2l} and
\ref{oracleconv} are too rough (non-optimal) when $M=M_n$ depends on
$n$ and the condition $M \le c_0 n$ is not satisfied. In the sequel,
we will apply those theorems with $M=M_n$ depending on $n$ and
satisfying $M_n/n \to 0$, as $n \to \infty$, so that the condition
$M \le c_0 n$ will obviously hold with any finite $c_0$ for $n$
large enough.

\begin{TH1}
\label{t1l} Let the integers $M\ge 2$ and $n\ge 1$ be such that $M
\le c_0 n$ where $c_0$ is a positive constant. Then there exist
probability densities $p_j \in L_2(\R^d),\, j=1,\dots,M,$ such that
for all estimators $T_n$ of $p$ we have
\begin{equation}
\label{eqt1l} \inf_{T_n} \sup_{p\in \mathcal{P}_0} \Big[R_n(T_n,p)-
\inf_{\lambda\in \R^M} \|\pp_{\lambda} - p\|^2\Big] \ge cLM/n
\end{equation}
where $c>0$ is a constant depending only on $c_0$.
\end{TH1}
{\sc Proof.} Set $r=M-1 \ge 1$ and fix $0<a<1$. Consider the
function $\tilde g$ defined for any $t \in \R$ by
$$
\tilde g(t)\triangleq \frac{aL}{2} \1_{\left[0,
    \frac{1}{Lr}\right]}(t)-\frac{aL}{2}\1_{\left(
    \frac{1}{Lr},\frac{2}{Lr}\right)}(t),
$$
where $\1_{A}(\cdot)$ denotes the indicator function of a set $A$.
Let $\{\tilde g_j\}_{j=1}^{r}$ be the family of functions defined by
$\tilde g_j(t)=\tilde g(t-2(j-1)/Lr),\ 1 \le j \le r$. Define also
the density $\tilde f(t)=(L/2)\1_{[0,2/L]}(t), \ t \in \R$. For
$x=(x_1,\ldots, x_d) \in \R^d$ consider the functions
$$
f(x)=\tilde f(x_1)\prod_{k=2}^d \1_{[0,1]}(x_k)\quad g_j(x)=\tilde
g_j(x_1)\prod_{k=2}^d \1_{[0,1]}(x_k), \ j=1,\ldots,r\,.
$$
Define the probability densities $p_j$ by $p_1=f$, $p_{j+1}=f+g_j,
\, j=1,\dots,M-1$.

Consider now the set of functions $\mathcal{Q}=
\{q_{\delta}:q_{\delta}=f + \sum_{j=1}^{r} \delta_jg_j, \,
\delta=(\delta_1,\ldots, \delta_r) \in  \{0,1\}^r\}$. Clearly, for
any $\delta \in \{0,1\}^r$, $q_{\delta}$ satisfies $\int_{\R^d}
q_{\delta}(x) \ud x =1$, $q_{\delta}\ge 0$ and
$\|q_{\delta}\|_\infty \le L$. Therefore $\mathcal{Q} \subset
\mathcal{P}_0$. Also, $\mathcal{Q} \subset \{\pp_\lambda, \lambda\in
\R^M\}$. Thus,
$$
 \inf_{T_n} \sup_{p\in \mathcal{P}_0} \Big[R_n(T_n,p)-
\inf_{\lambda\in \R^M} \|\pp_{\lambda} - p\|^2\Big] \ge \inf_{T_n} \sup_{p\in
\mathcal{Q}}R_n(T_n,p)\,.
$$
To prove that $ \inf_{T_n} \sup_{p\in \mathcal{Q}}R_n(T_n,p) \ge c
LM/n $ we check conditions \eqref{condlem1} of Lemma~\ref{lem1}. The
first condition in \eqref{condlem1} is obviously satisfied since
$$
\|g_j\|^2=\int_0^{\frac{2}{Lr}}{\tilde g}^2(t) \ud t
=\frac{a^2L}{2r}\, , \quad j=1,\dots,r.
$$
To check the second condition in \eqref{condlem1}, note that for
$j=1,\dots,r$ we have
\begin{equation*}
\begin{split}
h^2(f, f+g_j)
&=\frac{1}{2}\int_0^{\frac{2}{Lr}}\Big(\sqrt{L/2}-\sqrt{L/2+\tilde
g(t)}\,\Big)^2 \ud t \\
& =\frac{L}{4}\int_0^{\frac{2}{Lr}}\Big(1-\sqrt{1+(2/L)\tilde
g(t)}\,\Big)^2 \ud t\\
&=\frac{L}{4}\bigg[\frac{4}{Lr}-2
\int_0^{\frac{2}{Lr}}\sqrt{1+(2/L)\tilde
g(t)}\ud t \bigg]\\
&= \frac{1}{r} -\frac{1}{2r} \Big( \sqrt{1+a}+\sqrt{1-a} \Big)
\le \frac{a^2}{2r} \,\\
\end{split}
\end{equation*}
where we used the fact that $\sqrt{1+a}+\sqrt{1-a} \ge 2-a^2$ for
$|a|\le 1$.
Define now $\tilde c_0=\max(c_0,3)$ and choose $a^2=M/(\tilde c_0
n)\le 1$. Then $a^2/(2r)\le (\tilde c_0n)^{-1}$ for $M\ge 2$.
Applying Lemma~\ref{lem1} with $\beta=(\tilde c_0n)^{-1}$ and $\DS
\alpha=\frac{ML}{2\tilde c_0 nr}$ we get
$$
\inf_{T_n} \sup_{p \in {\cal C}}R_n(T_n,p) \ge \frac{1}{8\tilde c_0}
\bigg(1-\sqrt{\frac{2}{\tilde c_0}}\,\bigg) \frac{LM}{n}\,.
$$
 \epr

Theorems~\ref{t2l} and~\ref{t1l} imply the following result.
\begin{cor}
Let the integers $M\ge 2$ and $n\ge 1$ be such that $M \le c_0 n$
where $c_0$ is a positive constant. Then $\psi_n(M)=LM/n$ is optimal
rate of linear aggregation over $\mathcal{P}_0$ (i.e. the optimal
rate of aggregation for $H= \R^M$ over $\mathcal{P}_0$), and $\tilde
p_n^{\bf L}$ defined in \eqref{11} is rate optimal aggregate for
$\R^M$ over $\mathcal{P}_0$.
\end{cor}
Consider now a lower bound for convex aggregation. We analyze here
only the case $H=\Lambda^M$. Other examples of convex sets $H$ can
be treated similarly.

\begin{TH1}
\label{t1c} Let the integers $M\ge 2$ and $n\ge 1$ be are such that
$M \le c_0 n$. Then there exist functions $p_j \in L_2(\R^d),
j=1,\dots,M,$ such that for all estimators $T_n$ of $p$ we have
\begin{equation}
\label{eqt1c} \inf_{T_n} \sup_{p\in \mathcal{P}_0} \Big[R_n(T_n,p)-
\inf_{\lambda\in \Lambda^M} \|\pp_{\lambda} - p\|^2\Big] \ge cLM/n
\end{equation}
where $c>0$ is a constant depending only on $c_0$.
\end{TH1}
{\sc Proof.} Consider the same family of densities $\mathcal{Q}$ as
defined in the proof of Theorem~\ref{t1l}. We may rewrite it in the
form $\mathcal{Q}=
\{q_{\delta}:q_{\delta}=\lambda_1Mf+\sum_{j=1}^{r} \lambda_{j+1}M
\delta_jg_j, \, \delta=(\delta_1,\ldots, \delta_r) \in
 \{0,1\}^r\}$ where $\lambda_j = 1/M$, $j=1, \ldots M$. Define now
$p_1=Mf$, $p_{j+1}=M(f+g_j), \, j=1,\dots,M-1$. Since $
\sum_{j=1}^{M} \lambda_j = 1 $ we have $\mathcal{Q} \subset \left\{
\pp_{\lambda} , \lambda \in \Lambda^M \right\}$. The rest of the
proof is identical to that of Theorem~\ref{t1l}. \epr

Theorems~\ref{oracleconv} and~\ref{t1c} imply the following result.
\begin{cor}
Let the integers $M\ge 2$ and $n\ge 1$ be such that $M \le c_0n$.
Then $\psi_n(M)=LM/n$ is optimal rate of convex aggregation over
$\mathcal{P}_0$ (i.e. the optimal rate of aggregation for $H=
\Lambda^M$ over $\mathcal{P}_0$), and $\tilde p_n^{\bf C}$ is rate
optimal aggregate for $H=\Lambda^M$ over $\mathcal{P}_0$.
\end{cor}

Inspection of the proofs of Theorems~\ref{t1c} and~\ref{t1l} reveals
that the least favorable functions $p_j$ used in the lower bound for
linear aggregation are uniformly bounded by $L$, whereas this is not
the case for least favorable functions in convex aggregation. It can
be shown that, for convex aggregation of functions which are
uniformly bounded by $L$, an elbow appears in the optimal rates of
aggregation, with the bound \eqref{eqt1c} still remaining valid for
$M\le \sqrt{n}$. This issue will be treated in a forthcoming paper
of the first author.

\section{Sample splitting and averaged aggregates}
\setcounter{equation}{0} \label{sec4}

We now come back to the original problem discussed in the
introduction. Let $\X_1^m$ denote a subsample of $\X^n=(X_1, \ldots,
X_n)$ of size $m\le n$ (training sample). Take $m<n$ and construct
estimators $\hat p_{m,1}, \ldots, \hat p_{m,M}$ of $p$ based on
$\X_1^m$. Then aggregate these estimators using the validation
subsample $\X_2^{\ell}$ of $\X^n$ of size $\ell=n-m$,
$$
\left(  \X_1^m, \X_2^{\ell} \right)=\X^n=(X_1, \ldots, X_n).
$$
For given $m<n$ the two subsamples can be obtained by different
splits. The choice of split is arbitrary, and it may influence the
result of estimation. In order to avoid the arbitrariness, we will
use a jackknife type procedure averaging the aggregates over
different splits. Define a \emph{split} $\mathcal{S}$ of the initial
sample $\X^n$ as a mapping
$$
\mathcal{S}\ : \ \X^n \mapsto \left(  \X_1^m, \X_2^{\ell} \right).
$$
Denote by $\X_{1,\mathcal{S}}^m, \X_{2,\mathcal{S}}^{\ell}$
subsamples obtained for a fixed split $\mathcal{S}$ and consider an
arbitrary set of splits $\mathbb{S}$. It can be, for example, the
set of all splits. Define $\tilde p^{\mathcal{S}}_n$ as a linear or
convex aggregate (${\tilde p}_n^{\bf L}$ or ${\tilde p}_n^{\bf C}$
respectively) based on the validation sample
$\X_{2,\mathcal{S}}^{\ell}$ and on the initial set of estimators
$p_j=\hat p^{\mathcal{S}}_{m,j}, j=1,\dots, M$, where each of $\hat
p^{\mathcal{S}}_{m,j}$'s is constructed from the training sample
$\X_{1,\mathcal{S}}^m$. Introduce the following {\it averaged
aggregate} estimator:
\begin{equation}
\label{defjack} \tilde p^{\mathbb{S}}_n \triangleq \frac{1}{{\rm
card}( \mathbb{S})} \sum_{\mathcal{S} \in \mathbb{S}} \tilde
p^{\mathcal{S}}_n.
\end{equation}
Let $H$ be either $\R^M$ or a convex compact subset of $\R^M$.
Define
$$
\Delta_{\ell,M}=\left\{
\begin{array}{cl}
  LM/\ell & \mbox{if} \ H=\R^M,\\
  4LM/\ell & \mbox{if} \ H \ \mbox{is a convex compact subset of $\R^M$}.\\
\end{array}
\right.
$$
We get the following corollary of Theorems~\ref{t2l} and
\ref{oracleconv}.

\begin{cor}
\label{c4.1} Let $m<n$, $\ell=n-m$, and let $H$ be either $\R^M$ or
a convex compact subset of $\R^M$. Let $\mathbb{S}$ be an arbitrary
set of splits. Assume that $\hat p^{\mathcal{S}}_{m,1}, \ldots, \hat
p^{\mathcal{S}}_{m,M}\in L_2(\R^d) $ for fixed
$\X_{1,\mathcal{S}}^m$, $\forall \mathcal{S}\in\mathbb{S}$, and that
$p \in \mathcal{P}_0$. Then the averaged aggregate~\eqref{defjack}
satisfies
\begin{equation}
\label{jack} R_n(\tilde p^{\mathbb{S}}_n ,p)
 \leq\inf_{\lambda \in H} R_m \Big( \sum_{j=1}^M \lambda_j \hat
p_{m,j} ,p \Big) + \Delta_{\ell,M}
\end{equation}
for any integers $M\ge 2$ and $n\ge 1$.
\end{cor}
{\sc Proof.} For any fixed $\mathcal{S}\in \mathbb{S}$ and for a
fixed training subsample $\X_{1,\mathcal{S}}^m$ inequalities
\eqref{bslin} and \eqref{bb} imply
\begin{equation}
\label{unifnot} E^{\ell,\mathcal{S}}_p\|\tilde
p^{\mathcal{S}}_n-p\|^2 \le \min_{\lambda\in H}\Big\|\sum_{j=1}^M
\lambda_j \hat p_{m,j} - p\Big\|^2 + \Delta_{\ell,M}, \quad \forall
\ p \in \mathcal{P}_0,
\end{equation}
where $E^{\ell,\mathcal{S}}_p$ denotes the expectation w.r.t. the
distribution of the validation sample $\X_{2,\mathcal{S}}^{\ell}$
when the true density is $p$. Taking expectations of both sides of
\eqref{unifnot} w.r.t. the training sample $\X_{1,\mathcal{S}}^m$ we
get
\begin{equation}
\label{unifnot1} R_n(\tilde p^{\mathcal{S}}_n ,p)
 \leq\inf_{\lambda \in H} R_m \Big( \sum_{j=1}^M \lambda_j \hat
p_{m,j} ,p \Big) + \Delta_{\ell,M}.
\end{equation}
The right hand side here does not depend on $\mathcal{S}$. By
Jensen's inequality,
$$
R_n(\tilde p^{\mathbb{S}}_n ,p)
 \leq
 \frac{1}{{\rm card}( \mathbb{S})}
\sum_{\mathcal{S} \in \mathbb{S}} R_n(\tilde p^{\mathcal{S}}_n ,p).
$$
This and \eqref{unifnot1} yield \eqref{jack}. \epr

\section{Kernel aggregates for density estimation}
\setcounter{equation}{0} \label{sec5}

Here we apply the results of the previous sections to aggregation of
kernel density estimators. Let $\hat p_{m,h}$ denote a kernel
density estimator based on $\X_1^m$ with $m\le n$,
\begin{equation}
\label{defkernel}
\hat p_{m,h}(x) \triangleq \frac{1}{mh^d} \sum_{i=1}^n K \left(\frac{X_i
    -x}{h} \right)\1_{\{\X_1^m\}}(X_i), \quad  x \in \R^d,
\end{equation}
where $h>0$ is a bandwidth and $K \in L_2(\R^d)$ is a kernel. The
notation $\hat p_{m,h}$ is slightly inconsistent with $\hat p_{m,j}$
used above but this will
not cause ambiguity in what follows. In order to cover such
examples as the sinc kernel we will not assume that $K$ is
integrable.

Define $h_0 =(n \log n)^{-1/d}$, $a_n= a_0/ \log n$, where $a_0>0$
is a constant, and $M$ such that
$$
M-2=\max\left\{ j \in \N:h_0(1+a_n)^j<1\right\}.
$$
It is easy to see that $M \leq c_4 (\log n)^2$, where $c_4>0$ is a
constant depending only on $a_0$ and $d$. Consider a grid
$\mathcal{H}$ on $[0,1]$ with a weakly geometrically increasing
step:
$$
\mathcal{H} \triangleq \left\{h_0, h_1, \ldots, h_{M-1}  \right\}\,,
$$
where $h_j=(1+a_n)^jh_0, \, j=1,\dots, M-2,$ and $h_{M-1}=1$. Fix
now an arbitrary family of splits $\mathbb{S}$ such that, for $n \ge
3$,
$$
m = \lfloor n \left(1 - (\log n)^{-1} \right) \rfloor \quad {\rm
and} \quad \ell=n-m \ge \frac{n}{\log n} \, ,
$$
where $\lfloor x \rfloor$ denotes the integer part of $x$.

Define $\tilde p^{\mathbb{S},K}_n$ as the linear or convex (with
$H=\Lambda^M$) averaged aggregate $\tilde p^{\mathbb{S}}_n$ where
the initial estimators are taken in the form $p_j = \hat
p_{m,h_{j-1}}, j=1,\dots, M,$ with $\hat p_{m,h}$ given by
\eqref{defkernel}. Since $\Delta_{\ell,
  M} \le 4LM/\ell$ we get from \eqref{jack} that, under the
  assumptions of Corollary~\ref{c4.1},
\begin{equation}
\label{grid} R_n(\tilde p^{\mathbb{S},K}_n ,p) \leq \min_{h \in
\mathcal{H}}R_m(\hat p_{m,h} ,p)+ \Delta_{\ell,M} \leq \min_{h \in
\mathcal{H}}R_m(\hat p_{m,h} ,p) + \frac{4c_4(\log n)^3}{n} \, .
\end{equation}
We now give a theorem that extends~\eqref{grid} to the $n$-sample oracle risk
$\inf_{h>0}R_n(\hat p_{n,h} ,p)$ instead of $\min_{h \in
\mathcal{H}}R_m(\hat p_{m,h} ,p)$. Denote by $\mathcal{F}[f]$ the
Fourier transform defined for $f \in L_2(\R^d)$ and normalized in
such a way that its restriction to $f \in L_2(\R^d)\cap L_1(\R^d)$
has the form $\mathcal{F}[f](t)= \int_{\R^d} \e^{i
  x^T t} f(x) \ud x, \, t\in \R^d$.
In the sequel $\varphi=\mathcal{F}[p]$ denotes the characteristic
function associated to $p$.

\begin{TH1}
\label{kernoracle}
Assume that $p$ satisfies $\|p\|_{\infty} \leq L$ with $0<L<\infty$ and let
 $K \in L_2(\R^d)$ be a kernel such that a version of its Fourier transform
 $\mathcal{F}[K]$ takes values in $[0,1]$ and satisfies the monotonicity condition
 $\mathcal{F}[K](h't)\ge
 \mathcal{F}[K](ht), \ \forall \, t \in \R^d,\, h>h'>0$.
Then there exists an integer $n_0=n_0(L,\|K\|)\ge 4$ such that for
$n \ge n_0$ the averaged aggregate $\tilde p^{\mathbb{S},K}_n$
satisfies the oracle inequality
\begin{equation}
\label{oraclekernel} R_n(\tilde p^{\mathbb{S},K}_n ,p) \leq
\left(1+c_{5}(\log n)^{-1}\right)\inf_{h >0}R_n(\hat p_{n,h} ,p)+
c_{6}\frac{ (\log n)^3}{n} \, ,
\end{equation}
where $c_{5}$ is a positive constant depending only on $d$ and
$a_0$, and $c_{6}>0$ depends only on $L,\|K\|, d$ and $a_0$.
\end{TH1}
{\sc Proof.} Assume throughout that $n\ge 4$. First note that
\eqref{oraclekernel} deduces from \eqref{grid} and from the
following two inequalities that we are going to prove below:
\begin{equation}
\label{eqproof1} \inf_{h \in [h_0, h_{M-1}]}R_n(\hat p_{n,h},p) \le
\inf_{h>0} R_n(\hat p_{n,h},p)+\|K\|^2\frac{\log n}{n} \, ,
\end{equation}
\begin{equation}
\label{eqproof2} \min_{j=1, \ldots,M}R_m(\hat p_{m,h_{j-1}},p)\leq
\left(1 + c_{5}(\log
  n)^{-1}\right)  \inf_{h \in [h_0, h_{M-1}]}R_n(\hat p_{n,h},p)
  +\frac{c_{5}L}{ n \log n} \, .
\end{equation}
In turn, \eqref{eqproof1} follows if we show that
\begin{eqnarray}\label{i}
&& \inf_{h \in [h_0, h_{M-1}]}R_n(\hat p_{n,h},p) \le
\inf_{0<h<h_0}R_n(\hat p_{n,h},p),\\
\label{ii} && \inf_{h \in [h_0, h_{M-1}]}R_n(\hat p_{n,h},p) \le
\inf_{h
>
  h_{M-1}}R_n(\hat p_{n,h},p)+\|K\|^2\frac{\log n}{n} \, .
\end{eqnarray}
Thus, it remains to prove \eqref{eqproof2} -- \eqref{ii}. We will
use the following Fourier representation for MISE of kernel
estimators that can be easily obtained from Plancherel's formula (it
is a multivariate extension of the representation for $d=1$ given,
e.g., in Golubev (1992) and in Wand and Jones (1995), p.55):
\begin{equation}
\label{plancherel}
\begin{split}
R_n(\hat
p_{n,h} ,p) & =\frac{1}{(2 \pi)^d} \int_{\R^d} \Big( | 1-\mathcal{F}[K](ht) |^2 |
 \varphi (t ) |^2 \\
& \phantom{\frac{1}{(2 \pi)^d} \int_{\R^d} \Big( | 1-} +
 \frac{1}{n} \left(1-|\varphi (t ) |^2 \right)\big| \mathcal{F}[K](h
 t)\big|^2\Big) \ud t.
\end{split}
\end{equation}
Furthermore, using Plancherel's formula we get
\begin{equation}
\begin{split}
\label{plancherel1}
&\int_{\R^d} |\varphi(t)|^2 \ud t = (2 \pi)^d
\int_{\R^d}p^2(x) \ud x \le (2 \pi)^d L,\\
& \frac{1}{(2 \pi)^d}
\int_{\R^d}  |\mathcal{F}[K](ht) |^2 \ud t  = h^{-d}\|K\|^2, \
\forall \, h>0.
\end{split}
\end{equation}
{\bf Proof of \eqref{i}}. Using \eqref{plancherel},
\eqref{plancherel1} and the fact that $0 \le \mathcal{F}[K](t) \le
1, \, \forall t\in \R^d$, for any $h<h_0=(n \log n)^{-1/d}$ we
obtain
\begin{equation}
\label{i.1} R_n(\hat p_{n,h},p)  \ge \frac{1}{n(2 \pi)^d} \int_{\R^d}
\left(1-|\varphi (t ) |^2 \right)\big| \mathcal{F}[K](h
 t)\big|^2 \ud t \ge \frac{\|K\|^2}{n h^d}-\frac{L}{n}\ge
 \|K\|^2\log n-\frac{L}{n}.
\end{equation}
On the other hand, since $h_{M-1}=1$ we get
\begin{equation}
\label{i.2} R_n(\hat p_{n,h_{M-1}},p)  \le  \frac{1}{(2 \pi)^d}
\int_{\R^d} \left( | 1-\mathcal{F}[K](t) |^2 |
 \varphi (t ) |^2  +
 \frac{1}{n}\big| \mathcal{F}[K](
 t)\big|^2\right) \ud t \le  L + \frac{\|K\|^2}{n}.
\end{equation}
The right hand side of \eqref{i.1} is larger than that of
\eqref{i.2} for $n \ge n_0$, where $n_0$ depends only on $L$ and
$\|K\|$. Thus, \eqref{i} is valid for $n \ge n_0$.

\noindent {\bf Proof of \eqref{ii}}. Clearly, \eqref{ii} follows if
we show that
$$
R_n(\hat p_{n,h'},p) \le \inf_{h >
  h_{M-1}}R_n(\hat p_{n,h},p)+\|K\|^2\frac{\log n}{n}
$$
for $h'=(\log n)^{-1/d} \in [h_0, h_{M-1}]$. To prove this
inequality, first note that, by the monotonicity of
$h\mapsto\mathcal{F}[K](ht)$,  we have
$$
 \int_{\R^d}  | 1-\mathcal{F}[K](ht) |^2 |
 \varphi (t ) |^2\ud t
\ge \int_{\R^d}  | 1-\mathcal{F}[K](h't) |^2 |
 \varphi (t ) |^2 \ud t, \quad \forall h>h_{M-1}.
$$
This, together with \eqref{plancherel} and the second equality in
\eqref{plancherel1}, yields that, for any $h>h_{M-1}$,
\begin{equation*}
R_n(\hat p_{n,h},p) \ge R_n(\hat p_{n,h'},p)-\frac{1}{n(2 \pi)^d} \int_{\R^d} \left(1-|\varphi (t ) |^2 \right)\big| \mathcal{F}[K](h'
 t)\big|^2\ud t \ge  R_n(\hat p_{n,h'},p)-\|K\|^2\frac{\log n}{n} \, .
\end{equation*}
{\bf Proof of \eqref{eqproof2}}. We will show that for any $h\in
[h_0, h_{M-1}]$ one has
\begin{equation}
\label{pr2} R_m(\hat p_{m,\overline h} ,p) \leq \left(1+c_{5}(\log
n)^{-1}\right)R_n(\hat p_{n,h} ,p) + \frac{c_5 L}{ n \log n}
\end{equation}
where $ \overline h \triangleq \max\{h_j:h_j \le h\}. $ Clearly,
this implies \eqref{eqproof2}. To prove \eqref{pr2}, note that
 if $h_j \leq h < h_{j+1}$ we have $\overline
h = h_j$, $h/h_j \le 1+a_n = 1+a_0/\log n$. Therefore,
\eqref{plancherel} and the monotonicity of
$h\mapsto\mathcal{F}[K](ht)$ imply
\begin{equation*}
\begin{split}
R_m(\hat p_{m,h_j} ,p) & = \frac{1}{(2 \pi)^d}\int_{\R^d}\left(
\left[1 - \mathcal{F}[K](h_j t) \right]^2|\varphi(t)|^2
 + \frac{1}{m} \left[ \mathcal{F}[K](h_j t) \right]^2 \right)\ud t\\
& \phantom{\frac{1}{(2 \pi)^d}\int_{\R^d} \left[1 - -- \right]^2} -
\frac{1}{(2 \pi)^d m}\int_{\R^d}|\varphi(t)|^2
  \left[ \mathcal{F}[K](h_j t) \right]^2\ud t\\
& \leq \frac{1}{(2 \pi)^d }\int_{\R^d}\Big( \left[1 -
\mathcal{F}[K](h t) \right]^2|\varphi(t)|^2
 + \frac{1}{n} \left[ \mathcal{F}[K](h t) \right]^2 \frac{nh^d}{mh_j^d}
\Big)\ud t\\
& \phantom{\frac{1}{(2 \pi)^d}\int_{\R^d} \left[1 - -- \right]^2} -
\frac{1}{(2 \pi)^d n}\int_{\R^d}|\varphi(t)|^2
  \left[ \mathcal{F}[K](h t) \right]^2\ud t\\
& \leq \frac{nh^d}{mh_j^d}
  R_n(\hat p_{n,h} ,p)+ \left(\frac{nh^d}{mh_j^d} - 1\right)
 \frac{1}{(2 \pi)^d n}\int_{\R^d}|\varphi(t)|^2
  \left[ \mathcal{F}[K](h t) \right]^2\ud t.
\end{split}
\end{equation*}
Using here the fact that $(n/m)(h/h_j)^d\le (1-(\log
n)^{-1}-n^{-1})(1+a_0/\log n)^d\le 1+c_5(\log n)^{-1}$ for $n\ge 4$
and for a constant $c_5>0$ depending only on $d$, $a_0$, and
applying \eqref{plancherel1} we get \eqref{pr2}. \epr

\begin{cor}
\label{corkernoracle} Let the assumptions of Theorem~\ref{kernoracle} be satisfied, and let $\inf_{h>0}R_n(\hat p_{n,h}
,p) \ge cn^{-1+\alpha}$, for some $c>0,\alpha>0$. Then
\begin{equation}
\label{exactoracle} R_n(\tilde p^{\mathbb{S},K}_n ,p) \leq \inf_{h
>0}R_n(\hat p_{n,h} ,p)(1+o(1)),\qquad n\to \infty.
\end{equation}
\end{cor}
Using the argument as in Stone (1984) it is not hard to check that
the assumption of Corollary~\ref{corkernoracle} is valid for any
non-negative kernel. In the one-dimensional case it also holds for
any kernel satisfying the conditions of Lemma~4.1 in Rigollet
(2006). On the difference to Rigollet (2006),
Corollary~\ref{corkernoracle} applies to multidimensional density
estimation.

Theorem~\ref{kernoracle} and Corollary~\ref{corkernoracle} show that
linear or convex aggregate $\tilde p^{\mathbb{S},K}_n$ mimics the
best kernel estimator, without being itself in the class of kernel
estimators with data-driven bandwidth. Another method with such a
property has been suggested recently by Rigollet (2006) in the
one-dimensional case; it is based on a block Stein procedure in the
Fourier domain.

The results of this section can be compared to the work on
optimality of bandwidth selection in the $L_2$ sense for kernel
density estimation. A key reference is the theorem of Stone (1984)
establishing that, under some assumptions,
$$
\lim_{n \to \infty} \frac{\|\hat p_{n,h_n}-p\|^2}{\inf_{h>0}\|\hat
  p_{n,h}-p\|^2}=1,\quad \textrm{with probability 1},
$$
where $h_n$ is a data-dependent bandwidth chosen by
cross-validation. Our results are of a different type, because they
treat convergence of expected risk rather than almost sure
convergence. In addition, we provide oracle inequalities with
precisely defined remainder terms that hold under mild assumptions
on the density and on the kernel. Unlike Stone (1984), we do not
require the one-dimensional marginals of the density $p$ to be
uniformly bounded. Wegkamp (1999) considers model selection approach
to bandwidth choice for kernel density estimation. His main result
is of the form of
  \eqref{exactoracle} with a model selection kernel estimator
  in place of $\tilde p^{\mathbb{S},K}_n$, but it is valid for bounded,
  nonnegative, Lipschitz kernels with compact support
  (similar assumptions on $K$ are imposed by Stone (1984)). Our
result covers kernels with unbounded support, for example, the
Gaussian and Silverman's kernels that are often implemented, and
Pinsker's kernel that gives sharp minimax adaptive estimators on
Sobolev classes (cf. Section~\ref{sec6} below). In a recent work of Dalelane
(2004) the choice of bandwidth and of the kernel by
  cross-validation is investigated for the one-dimensional case ($d=1$).
  She provides an oracle
  inequality similar to \eqref{oraclekernel} with a remainder term
  of the order
  $n^{\delta-1},\ 0<\delta<1$, instead of $(\log n)^3/n$ that we have here.

 All these papers consider the model selection approach, i.e., they
study
  estimators with a single data-driven bandwidth chosen from a set of
  candidate bandwidths. Our approach is different since we estimate the
  density by a linear or convex combination of kernel estimators with bandwidths in the
  candidate set. Simulations (see Section~\ref{sec7} below) show that in most cases one of these estimators
  gets highly dominant weight in the resulting mixture. However, inclusion
  of other estimators with some smaller weights allows one to treat more
  efficiently densities with inhomogeneous smoothness.

\section{Sharp minimax adaptivity
of kernel aggregates}
  \setcounter{equation}{0} \label{sec6}

In this section we show that the kernel aggregate defined in Section~\ref{sec5} is sharp minimax adaptive over a scale of Sobolev classes
of densities.

For any $\beta>0,\ Q>0$ and any integer $d \ge 1$ define the Sobolev
classes of densities on $\R^d$ by
$$
\Theta(\beta,Q) \triangleq \left\{ p :\R^d\to \R \, \Big| \ p\ge 0,
\int_{\R^d}p(x)\ud x=1,
  \ \int_{\R^d}\|t\|^{2\beta}_d |\varphi(t)|^2
  \ud t \le Q \right\},
$$
where $\|\cdot\|_d$ denotes the Euclidean norm in $\R^d$ and
$\varphi=\mathcal{F}[p]$. Consider the Pinsker kernel $K_{\beta}$,
i.e. the kernel having the Fourier transform
$$
\mathcal{F}[K_{\beta}](t) \triangleq
\left(1-\|t\|^{\beta}_d\right)_+, \quad t \in \R^d,
$$
where $x_+=\max(x,0)$. Set
\begin{equation}
\label{cstar} C^*=\frac{[Q(2\beta +d)]^{\frac{d}{2\beta+d}}}
{d(2\pi)^d} \left(\frac{\beta
S_d}{\beta+d}\right)^{\frac{2\beta}{2\beta+d}}
\end{equation}
where $S_d=2\pi^{d/2}/\Gamma(d/2)$ is the surface of a sphere of
radius 1 in $\R^d$. For $d=1$ the value $C^*$ equals to the Pinsker
constant [Pinsker (1980), see also Tsybakov (2004), Chapter 3].

\begin{cor}
\label{corkernel} For any integer $d \ge 1$ and any $\beta>d/2, \
Q>0$, the averaged linear or convex kernel aggregate $\tilde
p^{\mathbb{S},K_{\hspace{-0.3mm} \beta}}_n$ defined in Section~\ref{sec5} satisfies
$$
\sup_{p \in \Theta(\beta,Q)}R_n(\tilde p^{\mathbb{S},K_{\hspace{-0.3mm} \beta}}_n
,p)\leq C^*n^{-\frac{2\beta}{2\beta +d}}(1+o(1)), \qquad n \to
\infty,
$$
where $C^*$ is defined in \eqref{cstar}.
\end{cor}
{\sc Proof.} Denote by $\hat p_{n,h}$ the kernel density estimator
defined in \eqref{defkernel} with $m=n$ and $K = K_{\beta}$.
 Using \eqref{plancherel} and the fact that $0 \le
\mathcal{F}[K_{\beta}](t) \le 1, \, \forall \, t \in \R^d$, we get
\begin{equation}
\label{pinsk1}
\begin{split}
R_n(\hat p_{n,h} ,p)& \le \frac{1}{(2 \pi)^d} \int_{\R^d} \left( |
1-\mathcal{F}[K_{\beta}](ht) |^2 |
 \varphi (t ) |^2  +
 \frac{1}{n} \big| \mathcal{F}[K_{\beta}](h
 t)\big|^2\right) \ud t \\
& \le \frac{1}{(2 \pi)^d}\left(Qh^{2 \beta}+
 \frac{1}{n} \int_{\R^d}  \big| \mathcal{F}[K_{\beta}](h
 t)\big|^2 \ud t\right), \forall \, h>0,\, p\in \Theta(\beta,Q).\\
\end{split}
\end{equation}
Now, choose $h$ satisfying
\begin{equation}
\label{balance} \int_{\R^d}
\|t\|^{\beta}_d\mathcal{F}[K_{\beta}](ht)\ud t=Qnh^{\beta}.
\end{equation}
The solution of \eqref{balance} is
$$
h=D^*n^{-\frac{1}{2\beta+d}} \quad \textrm{where} \quad
D^*=\left(\frac{\beta
S_d}{Q(\beta+d)(2\beta+d)}\right)^{\frac{1}{2\beta+d}}.
$$
With $h$
satisfying \eqref{balance}, inequality~\eqref{pinsk1} becomes
\begin{equation*}
\begin{split}
R_n(\hat p_{n,h} ,p)& \leq  \frac{1}{(2 \pi)^dn} \int_{\R^d}
\mathcal{F}[K_{\beta}](ht)\left[
  \mathcal{F}[K_{\beta}](ht)+\|ht\|^{\beta}_d\right] \ud t\\
& =  \frac{1}{(2 \pi)^dnh^d}\int_{\R^d}
\mathcal{F}[K_{\beta}](t) \ud t
\\&
=  \frac{1}{(2 \pi)^dnh^d}\int_0^1\left(1-r^{\beta}\right)r^{d-1}S_d
\ud r\\
&= C^*n^{-\frac{2\beta}{2\beta +d}}.
\end{split}
\end{equation*}
Thus,
\begin{equation}\label{6.4}
\inf_{h>0} R_n(\hat p_{n,h} ,p) \leq C^*n^{-\frac{2\beta}{2\beta
+d}},  \quad \forall \ p \in \Theta(\beta,Q).
\end{equation}
Note that the kernel $K=K_{\beta}$ satisfies the conditions of
Theorem~\ref{kernoracle}, and it is easy to see that for $\beta>d/2$
there exists a constant $0<L<\infty$ such that $\|p\|_{\infty} \leq
L$ for all $p \in \Theta(\beta,Q)$. Thus, \eqref{oraclekernel}
holds, and to prove the corollary it suffices to take suprema of
both sides of \eqref{oraclekernel} over $p \in \Theta(\beta,Q)$ and
to use \eqref{6.4}. \epr

Along with Corollary~\ref{corkernel}, for any $\beta>d/2, \ Q>0$ the
following lower bound holds:
\begin{equation}
\label{lowerSobol}
 \inf_{T_n}\sup_{p \in \Theta(\beta,Q)}R_n(T_n ,p)\geq
C^*n^{-\frac{2\beta}{2\beta +d}}(1+o(1)), \qquad n \to \infty,
\end{equation}
 where
$C^*$ is defined in \eqref{cstar} and $\inf_{T_n}$ denotes the
infimum over all estimators of $p$.
 For $d=1$ the bound \eqref{lowerSobol} can be deduced from the results of
Golubev (1991, 1992); it is also proven explicitly in Schipper
(1996) (for integer $\beta$) and in Rigollet (2006), Dalelane (2004)
(for all $\beta>1/2$). For $d>1$ the bound \eqref{lowerSobol} can be
found for a slightly different but essentially analogous minimax
setup in Efromovich (2000). 
Corollary~\ref{corkernel} and the lower
bound \eqref{lowerSobol} imply that the estimator $\tilde
p^{\mathbb{S},K_{\beta}}_n$ is asymptotically minimax in the exact
sense (with the constant) over the Sobolev class of densities
$\Theta(\beta,Q)$ and is adaptive to $Q$ for any given $\beta$.
However, $\tilde p^{\mathbb{S},K_{\beta}}_n$ is not adaptive to the
unknown smoothness $\beta$ since the Pinsker kernel $K_{\beta}$
depends on $\beta$.

To get adaptation to $\beta$, we need to push aggregation one step
forward: we will aggregate
 kernel density estimators not only for different bandwidths but also
 for different kernels. To this end, we refine the notation $\hat p_{n, h}$ of
\eqref{defkernel} to $\hat p_{n, h, K}$, indicating the dependence
of the density estimator both on kernel $K$ and bandwidth $h$.
For a family of $N\ge 2$ kernels, $\mathcal{K}=\{K_{(1)}, \ldots,
K_{(N)}\}$, define $\tilde p^{\mathbb{S},\mathcal{K}}_n$ as the
linear or convex averaged aggregate where the initial estimators are
taken in the collection of kernel density estimators $\{\hat
p_{n,h,K}, K \in \mathcal{K}, h \in \mathcal{H}\}$. Thus, we
aggregate now $NM$ estimators instead of $M$. The following
corollary is obtained by the same argument as Theorem
\ref{kernoracle}, by merely inserting the minimum over ${K \in
\mathcal{K}}$ in the oracle inequality and by replacing $\|K\|$ with
its upper or lower bounds in the remainder terms.

\begin{cor}
\label{c6.2} Assume that $p$ satisfies $\|p\|_{\infty} \le L$ with
$0 < L < \infty$ and let $\mathcal{K}=\{K_{(1)}, \ldots, K_{(N)}\}$
be a family of kernels satisfying the assumptions of Theorem
\ref{kernoracle} and such that there exist constants $0< \ubar c <
\bar c < \infty$ with $\ubar c < \|K_{(j)}\| < \bar c, \ j=1, \dots,
N$. Then there exists an integer $n_1=n_1(L,\ubar c, \bar c)\ge 4$
such that for $n \ge n_1$ the averaged aggregate $\tilde
p^{\mathbb{S},\mathcal{K}}_n$ satisfies the oracle inequality
\begin{equation}
\label{oraclekernel2} R_n(\tilde p^{\mathbb{S},\mathcal{K}}_n ,p) \leq
\left(1+c_{5}(\log n)^{-1}\right)\min_{K \in \mathcal{K}} \inf_{h >0}R_n(\hat p_{n,h,K} ,p)+
c_{7}\frac{ N(\log n)^3}{n} \, ,
\end{equation}
where $c_{5}>0$ is the same constant as in Theorem~\ref{kernoracle}, and $c_{7}>0$ depends only on
$L,\ubar c, \bar c, d$ and $a_0$.
\end{cor}

Consider now a particular family of kernels $\mathcal{K}$. Define
$\mathcal{B}=\{\beta_1, \ldots, \beta_N\}$ where $\beta_1=d/2$,
$\beta_j=\beta_{j-1}+N^{-1/2}, \, j=2, \ldots,N$, and let
$\mathcal{K_{\mathcal{B}}}=\{K_b,\ b \in \mathcal{B}\}$ be a family
of Pinsker kernels indexed by $b\in \mathcal{B}$. We will later
assume that $N=N_n\to\infty$, as $n\to\infty$, but for the moment
assume that $N\ge 2$ is fixed. Note that
$\mathcal{K}=\mathcal{K}_{\mathcal{B}}$ satisfies the assumptions of
Corollary~\ref{c6.2}. In fact,
\begin{equation*}\label{ksq}
\|K_\beta\|^2 = S_d Q_d(\beta) \quad \mbox{where} \quad  Q_d(
\beta)=\frac{1}{d}-\frac{2}{\beta+d}+\frac{1}{2\beta +d} \, ,
\end{equation*}
and
\begin{equation}\label{ksq1}
\frac{1}{6d}\le Q_d(\beta) \le \frac{1}{d} \, , \ \ \ \forall \
\beta\ge d/2.
\end{equation}
Thus, the oracle inequality \eqref{oraclekernel2} holds with
$\mathcal{K}=\mathcal{K}_{\mathcal{B}}$.
 We will now prove that, under the assumptions of Corollary \ref{c6.2} the
linear or convex aggregate $\tilde
p^{\mathbb{S},\mathcal{K}_{\mathcal{B}}}_n$ with the initial
estimators in $\{\hat p_{n,h,K}, K \in \mathcal{K}_{\mathcal{B}}, h
\in \mathcal{H}\}$ satisfies the following inequality where $\beta$
in the oracle risk varies continuously:
\begin{equation}
\label{oraclekernel3}
R_n(\tilde
p^{\mathbb{S},\mathcal{K}_{\mathcal{B}}}_n,p) \leq
\Big(1+\frac{c_{5}}{\log n}\Big) \Big(1+\frac{6}{\sqrt{N}}\Big)
\inf_{\substack{ h
>0\\d/2<\beta<\beta_N }}R_n(\hat p_{n,h,K_{\beta}} ,p)+ c_{8}\frac{
N(\log n)^3}{n}.
\end{equation}
Fix $\beta \in (d/2,\beta_N), Q>0$ and $p \in \Theta(\beta, Q)$.
Define $\bar \beta=\min\{\beta_j \in \mathcal{B}: \beta_j>\beta\}$. In view of \eqref{oraclekernel2} with
$\mathcal{K}=\mathcal{K}_{\mathcal{B}}$, to prove
\eqref{oraclekernel3} it is sufficient to show that for any $h>0$
one has
\begin{equation}
\label{oraclekernel4} R_n(\hat p_{n,h,K_{\bar \beta}},p)\leq (1+
6N^{-1/2}) \left(R_n(\hat p_{n,h,K_{\beta}} ,p)+\frac{L}{n}\right).
\end{equation}
Using \eqref{plancherel} and the inequality $\bar \beta > \beta$ we
get
\begin{equation}
R_n(\hat p_{n,h,K_{\bar \beta}} ,p) \le R_n(\hat p_{n,h,K_{ \beta}}
,p)+\cI(\bar \beta)-\cI(\beta) \label{planch-pinsker}
\end{equation}
where
\begin{equation*}
\cI(\beta)\triangleq \frac{1}{(2 \pi)^dn}\int_{\R^d}
(1-\|ht\|_d^{\beta})_+^2 \ud t=\frac{\|K_{\beta}\|^2}{(2 \pi)^dnh^d}
= \frac{S_d}{(2 \pi)^dnh^d}Q_d(\beta).
\end{equation*}
Now, \mbox{$ Q_d(\bar \beta)=Q_d(\beta)+(\bar \beta -
\beta)Q_d'(b_0) $} for some $b_0\in [\beta,\bar \beta]$. Using
\eqref{ksq1} and the inequality $|Q_d'(\beta)|\le 1/d^2$ valid for
all $\beta>d/2$, we find that
$$Q_d(\bar \beta) \le Q_d(\beta)
+ 6(\bar \beta
  -\beta)Q_d(\beta)\le (1+
6N^{-1/2})Q_d(\beta)\,.
$$
Therefore,
\begin{equation}\label{planch1}
\cI(\bar \beta)\le (1+ 6N^{-1/2}) \cI(\beta).
\end{equation}
Also, in view of \eqref{plancherel} and \eqref{plancherel1} we have
\begin{equation}\label{planch2}
\cI(\beta)\le R_n(\hat p_{n,h,K_{\beta}} ,p)+\frac{L}{n}.
\end{equation}
Combining \eqref{planch-pinsker}, \eqref{planch1} and
\eqref{planch2} we obtain \eqref{oraclekernel4}, thus proving
\eqref{oraclekernel3}.

\begin{cor}
\label{c6.3} Assume that ${\rm Card}
(\mathcal{K}_{\mathcal{B}})=N_n$ where $\lim_{n\to\infty}N_n =
\infty$ and $\limsup_{n\to\infty} N_n /(\log n)^{\nu} < \infty$ for
some $\nu>0$. Then for any integer $d \ge 1$ and any $\beta>d/2, \
Q>0$, the averaged linear or convex kernel aggregate $\tilde
p^{\mathbb{S},\mathcal{K}_{\mathcal{B}}}_n$ satisfies
$$
\sup_{p \in \Theta(\beta,Q)}R_n(\tilde
p^{\mathbb{S},\mathcal{K}_{\mathcal{B}}}_n ,p)\leq
C^*n^{-\frac{2\beta}{2\beta +d}}(1+o(1)), \qquad n \to \infty,
$$
where $C^*$ is defined in \eqref{cstar}.
\end{cor}
{\sc Proof.} Fix $\beta >d/2, Q>0$. Let $n$ be large enough to
guarantee that $\beta< \beta_{N_n}$. Then the infimum on the right
in \eqref{oraclekernel3} is smaller or equal to
$C^*n^{-\frac{2\beta}{2\beta +d}}$ for all $p \in \Theta(\beta,Q)$
[cf. \eqref{6.4}]. To conclude the proof, it suffices to take
suprema of both sides of \eqref{oraclekernel3} over $p \in
\Theta(\beta,Q)$ and then pass to the limit as $n\to\infty$.\epr

Corollary~\ref{c6.3} and the lower bound \eqref{lowerSobol} imply
that the aggregate $\tilde
p^{\mathbb{S},\mathcal{K}_{\mathcal{B}}}_n$ is asymptotically
minimax in the exact sense (with the constant) over all Sobolev
classes of densities with $\beta>d/2$, $Q>0$, and thus it is sharp
adaptive (recall that its construction does not depend on the
parameters $Q$ and $\beta$ of the class).

\section{Simulations}
\label{sec7} \setcounter{figure}{0}

\begin{figure}[h]\centering
\begin{minipage}[t]{0.49\textwidth}
\centering
\includegraphics[width=\textwidth]{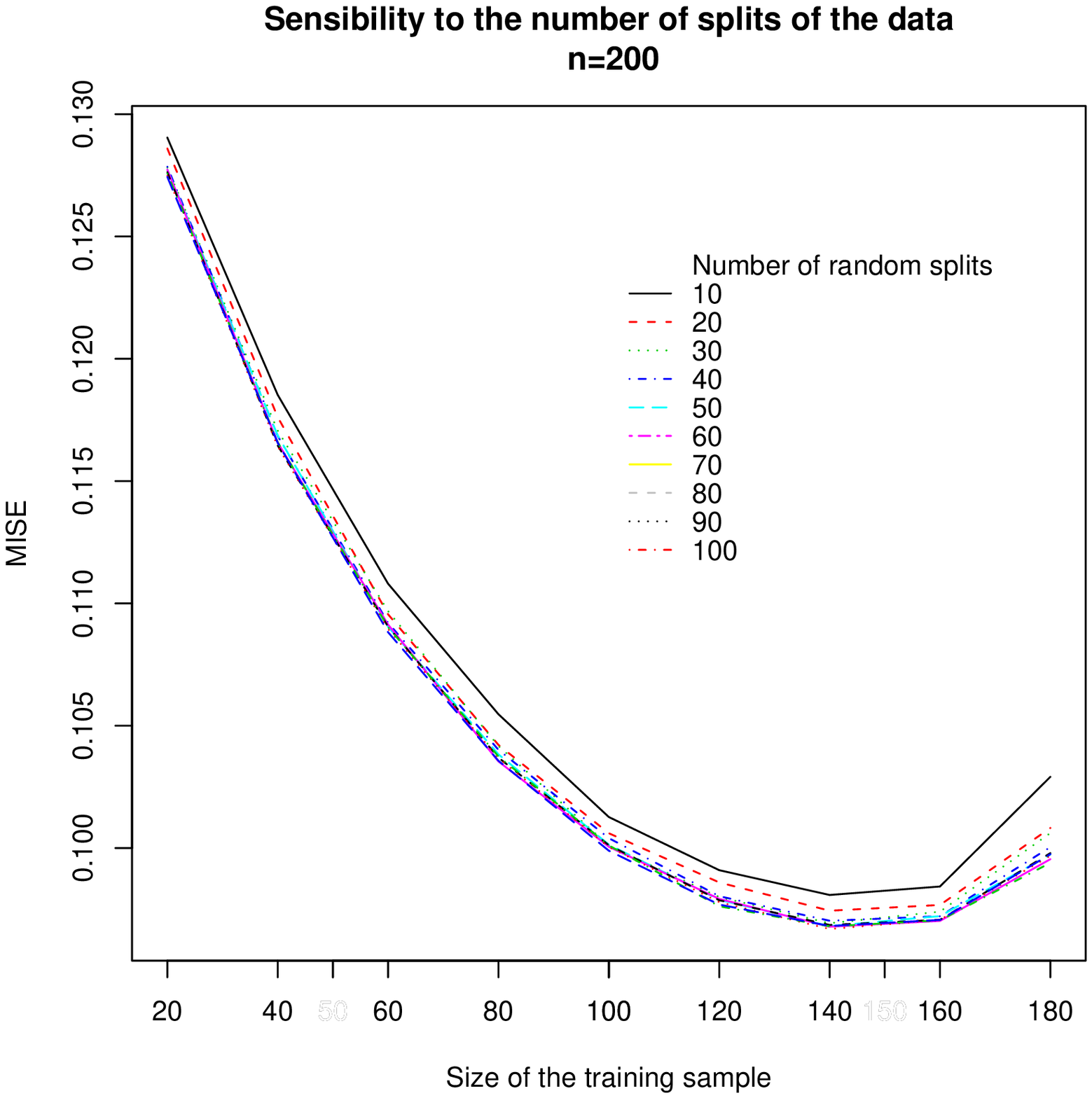}
\end{minipage}
\begin{minipage}[t]{0.49\textwidth}
\centering
\includegraphics[width=\textwidth]{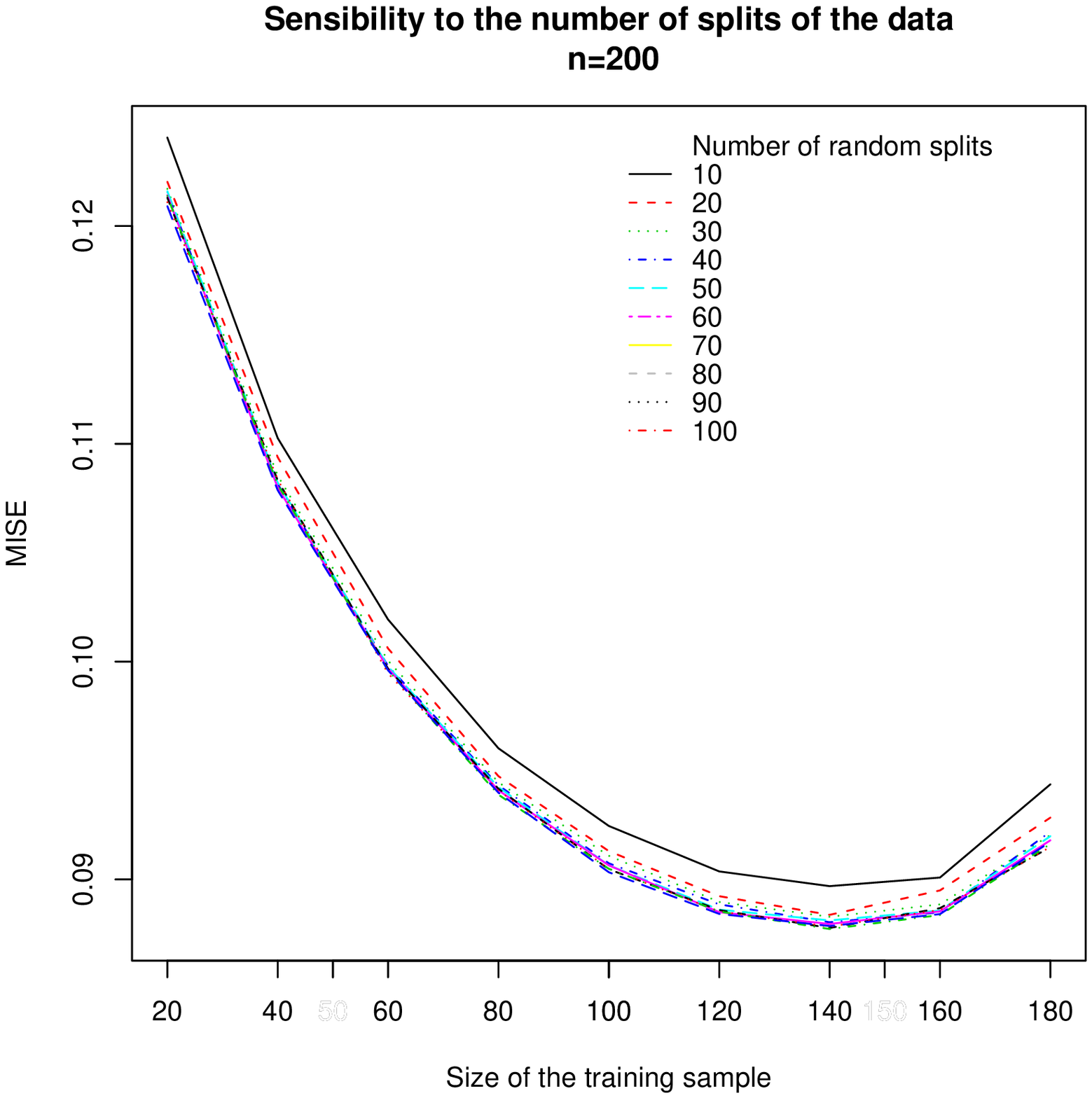}
\end{minipage}
\caption{Sensibility to the number of splits for dens1 (left) and
dens2 (right) }\label{sensplit}
\end{figure}

Here we discuss the results of simulations for the averaged convex
kernel aggregate with $H=\Lambda^M$ in the one-dimensional case. We
focus on convex aggregation because simulations of linear aggregates
show less numerical stability. The set of splits $\mathbb{S}$ is
reduced to 10 random splits of the sample since we observed that the
estimator is already stable for this number (cf.
Figure~\ref{sensplit}). In the default simulations each sample is
divided into two subsamples of equal sizes. The samples are drawn
from 6 densities that can be classified in the following three
groups.
\begin{itemize}
\item Common reference densities: the standard Gaussian density and the
  standard exponential density.
\item Gaussian mixtures from Marron and Wand (1992) that are
known to be difficult to estimate. We consider the Claw density and
the Smooth Comb density.
\item Densities with highly inhomogeneous smoothness. We consider
two densities referenced to as dens1 and dens2 that are both
mixtures of the standard Gaussian density $\varphi(\cdot)$ and of an
oscillating density. They are defined as
$$
0.5\varphi(\cdot)+0.5\sum_{i=1}^T \1_{\big(\frac{ 2(i-1)}{T}
\,,\frac{2i-1}{T}\big]}(\cdot)\,,
$$
where $T=14$ for dens1 and $T=10$ for dens2.

\end{itemize}
\noindent We used the procedure defined in Section~\ref{sec5} to
aggregate 6 kernel density estimators constructed with the Gaussian
$\n(0,1)$ kernel $K$ and with bandwidths $h$ from the set
$\mathcal{H} = \{0.001, 0.005, 0.01, 0.05, 0.1, 0.5\}$. This
procedure is further called \emph{pure kernel aggregation} and
quoted as {\tt AggPure}. Another estimator that we analyze is {\tt
AggStein} procedure: it aggregates 7 estimators, namely the same 6
kernel estimators as for {\tt AggPure} to which we add the block
Stein density estimator described in Rigollet (2006).  The
optimization problem \eqref{deflamconv} that provides aggregates is
solved numerically by a quadratic programming solver under linear
constraints: here we used the package \texttt{quadprog} of R. Our
simulation study shows that {\tt AggPure} and {\tt AggStein} have a
good performance for moderate sample sizes and are reasonable
competitors to kernel density estimators with common bandwidth
selectors.

We start the simulation by a comparison of the Monte-Carlo mean
integrated squared squared error (MISE) of {\tt AggPure} and {\tt
AggStein} with benchmarks. The MISE has been computed by averaging
integrated squared errors of 200 aggregate estimators calculated
from different samples of size 50, 100, 200 and 500. We compared the
performance of the convex aggregates and kernel estimators with
common data-driven bandwidth selectors and Gaussian $\n(0,1)$
kernel. The following bandwidth selectors are taken from the default
package {\tt stats} of the R software.
\begin{itemize}
\item DPI that implements the direct plug-in method of Sheather and Jones (1991) to
     select the bandwidth using pilot estimation of derivatives.
\item UCV and BCV that implement unbiased and biased
     cross-validation respectively (see, e.g., Wand and Jones (1995)).
\item Nrd0 that implements Silverman's rule-of-thumb [cf. Silverman (1986), page
     48]. It defaults the choice of bandwidth to 0.9 times the
     minimum of the standard deviation and the interquartile range
     divided by 1.34 times the sample size to the negative one-fifth
     power.
\end{itemize}
These descriptions correspond to the function {\tt
  bandwidth} in R which also allows for another choice of
  rule-of-thumb called {\tt Nrd}. It is a modification of {\tt Nrd0} given by Scott
  (1992), using factor 1.06 instead of 0.9. In our case, on the tested densities and sample
  sizes, this always leads to
  a MISE greater than that of {\tt Nrd0} except for the Gaussian density for which it is
  tailored. For this density, the performance of {\tt Nrd} is presented
  instead of that of {\tt Nrd0}.

The results are reported in Tables~\ref{miseGauExp} to
\ref{misedens777} where we included also the MISE of the block Stein
density estimator described in Rigollet (2006) and the oracle risk
which is defined as the minimum MISE of kernel density estimators
over the grid $\mathcal{H}$. It is, in general, greater than the
convex oracle risk, that is why it sometimes slightly exceeds the
MISE of convex aggregates or of other estimators that mimic more
powerful oracles for specific densities (such as DPI or {\tt Nrd}
for the Gaussian density).

\begin{table}[h]
{\small
\begin{minipage}[t]{0.49\textwidth}
\begin{tabular}[b]{|r|c|c|c|c|c|}
\hline
& 50 & 100 & 150 & 200 & 500 \\
\hline \hline
AggPure & 0.020 & 0.011 & 0.008 & 0.006 & 0.002\\
\hline
AggStein & 0.017 & 0.009 & 0.006 & 0.005 & 0.002 \\
\hline
Stein &  0.016 & 0.010 & 0.006 & 0.005 & 0.003\\
\hline
DPI &  0.011 & 0.006 & 0.005 & 0.004 & 0.002\\
\hline
UCV &  0.015 & 0.008 & 0.006 & 0.005 & 0.002\\
\hline
BCV &  0.009 & 0.006 & 0.004 & 0.003 & 0.002\\
\hline
Nrd &  0.010 & 0.006 & 0.004 & 0.003 & 0.002\\
\hline
Oracle & 0.008 & 0.005 & 0.004 & 0.004 & 0.003\\
\hline
\end{tabular}
\end{minipage}
\begin{minipage}[t]{0.49\textwidth}
\begin{tabular}[b]{|c|c|c|c|c|}
\hline
 50 & 100 & 150 & 200 & 500 \\
\hline \hline
 0.084 & 0.057 & 0.046 & 0.039 & 0.025 \\
\hline
 0.085 & 0.057 & 0.045 & 0.039 & 0.025 \\
\hline
 0.073 & 0.056 & 0.046 & 0.041 & 0.027 \\
\hline
 0.075 & 0.060 & 0.052 & 0.045 & 0.033 \\
\hline
 0.072 & 0.052 & 0.042 & 0.038 & 0.023 \\
\hline
 0.108 & 0.083 & 0.070 & 0.058 & 0.036 \\
\hline
 0.085 & 0.072 & 0.067 & 0.061 & 0.051 \\
\hline
 0.067 & 0.047 & 0.039 & 0.035 & 0.022 \\
\hline
\end{tabular}
\end{minipage}
\caption{\textit{MISE for the Gaussian (left) and the exponential
(right)
 densities} \label{miseGauExp}}
}
\end{table}

\begin{table}[h]
{\small
\begin{minipage}[t]{0.49\textwidth}
\begin{tabular}[b]{|r|c|c|c|c|c|}
\hline
& 50 & 100 & 150 & 200 & 500 \\
\hline \hline
AggPure & 0.058 & 0.041 & 0.034 & 0.029 & 0.014 \\
\hline
AggStein & 0.056 & 0.041 & 0.032 & 0.025 & 0.010 \\
\hline
Stein &  0.061 & 0.035 & 0.024 & 0.018 & 0.009\\
\hline
DPI &  0.059 & 0.052 & 0.050 & 0.048 & 0.043 \\
\hline
UCV & 0.063 & 0.043 & 0.032 & 0.026 & 0.012 \\
\hline
BCV & 0.058 & 0.052 & 0.051 & 0.050 & 0.046 \\
\hline
Nrd0 &  0.058 & 0.051 & 0.050 & 0.048 & 0.043\\
\hline
Oracle &  0.058 & 0.037 & 0.029 & 0.025 & 0.012\\
\hline
\end{tabular}
\end{minipage}
\begin{minipage}[t]{0.49\textwidth}
\begin{tabular}[b]{|c|c|c|c|c|}
\hline
 50 & 100 & 150 & 200 & 500 \\
\hline \hline
 0.064 & 0.042 & 0.034 & 0.029 & 0.017\\
\hline
 0.061 & 0.042 & 0.033 & 0.028 & 0.017 \\
\hline
 0.057 & 0.041 & 0.033 & 0.028 & 0.017\\
\hline
  0.070  & 0.054 & 0.046 & 0.042 & 0.029\\
\hline
  0.057 & 0.038 & 0.031 & 0.026 & 0.016\\
\hline
 0.101 & 0.083 & 0.066 & 0.055 & 0.027\\
\hline
 0.088 & 0.078 & 0.072 & 0.069 & 0.057\\
\hline
 0.064 & 0.038 & 0.030 & 0.025 & 0.016\\
\hline
\end{tabular}
\end{minipage}
\caption{\textit{MISE for the claw (left) and the smooth comb
(right)
 densities} \label{misecsc}}
}
\end{table}

\begin{table}[h]
{\small
\begin{minipage}[t]{0.49\textwidth}
\begin{tabular}[b]{|r|c|c|c|c|c|}
\hline
& 50 & 100 & 150 & 200 & 500 \\
\hline \hline
AggPure & 0.145 & 0.125 & 0.111 & 0.100 & 0.067\\
\hline
AggStein &  0.148 & 0.124  & 0.112 & 0.102 & 0.067\\
\hline
Stein &  0.152 & 0.143 & 0.140 & 0.138  & 0.132\\
\hline
DPI &  0.149 & 0.142 & 0.139 & 0.137 & 0.132 \\
\hline
UCV &  0.153 & 0.148 & 0.140 & 0.136 & 0.116 \\
\hline
BCV &  0.149 & 0.143 & 0.140 & 0.139 & 0.134 \\
\hline
Nrd0 &  0.149 & 0.141 & 0.138 & 0.137 & 0.133\\
\hline
Oracle & 0.148 & 0.144 & 0.142 & 0.133 & 0.067\\
\hline
\end{tabular}
\end{minipage}
\begin{minipage}[t]{0.49\textwidth}
\begin{tabular}[b]{|r|c|c|c|c|c|}
\hline
 50 & 100 & 150 & 200 & 500 \\
\hline \hline
 0.142 & 0.119 & 0.102 & 0.093 & 0.061 \\
\hline
  0.148 & 0.141 & 0.103 & 0.092 & 0.060\\
\hline
  0.154 & 0.143 & 0.140 & 0.137 & 0.132\\
\hline
  0.147 & 0.140 & 0.138 & 0.136 & 0.132 \\
\hline
  0.154 & 0.142 & 0.133 & 0.126 & 0.074\\
\hline
  0.146 & 0.141 & 0.139 & 0.138 & 0.134 \\
\hline
 0.146 & 0.140 & 0.137 & 0.136 & 0.132 \\
\hline
 0.145 & 0.128 & 0.109 & 0.101 & 0.062 \\
\hline
\end{tabular}
\end{minipage}
\caption{\textit{MISE for dens1 (left) and dens2 (right)}
\label{misedens777}} }
\end{table}

It is well known (see, e.g., Wand and Jones (1995)) that
  bandwidth selection by cross-validation (UCV) is unstable and leads
  too often to
  undersmoothing.
  The DPI and BCV
  methods were proposed in order to bypass the problem of
  undersmoothing. However, sometimes they
  lead to oversmoothing as in the case of the Claw density while convex
  aggregation works well. For the normal density DPI, BCV and Nrd are better, which comes
  as no surprise since these estimators are designed to estimate this density
  well. For the other densities that are more difficult to estimate these data
  driven bandwidth selectors do not provide good estimators whereas the
  aggregation procedures remain stable. The block Stein estimator
  performs well in all the cases except for the highly inhomogeneous
  densities (cf. Table 3).
In conclusion, the estimators {\tt AggPure} and {\tt AggStein} are
very robust, as compared to other tested procedures: they are not
far from the best performance for the densities that are easy to
estimate and they are clear winners for densities with inhomogeneous
smoothness for which other procedures fail.

\begin{figure}[h]\centering
\begin{minipage}[t]{0.49\textwidth}
\centering
\includegraphics[width=\textwidth]{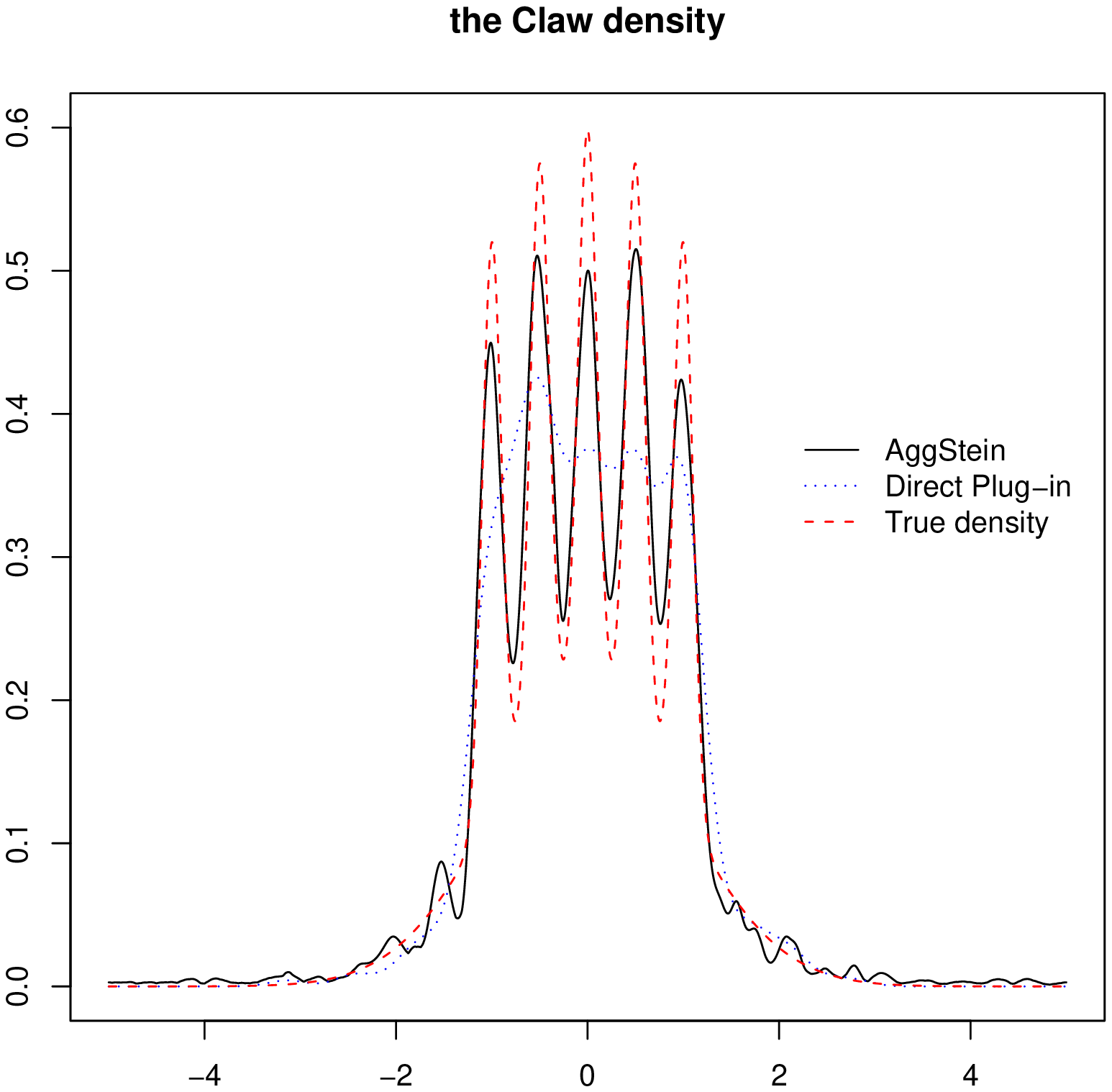}
\end{minipage}
\begin{minipage}[t]{0.49\textwidth}
\centering
\includegraphics[width=\textwidth]{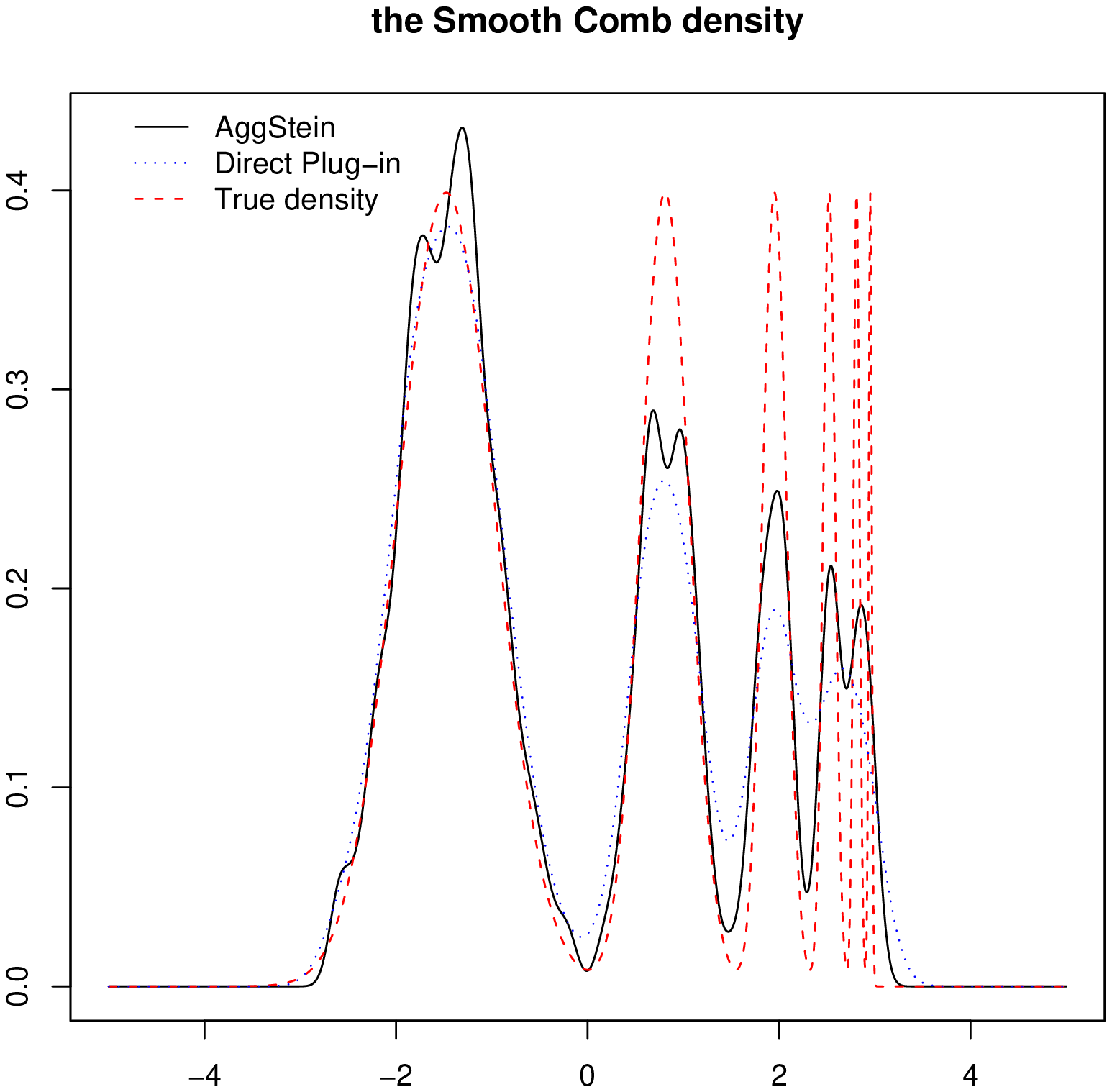}
\end{minipage}
\caption{The Claw and Smooth Comb densities }\label{compag}
\end{figure}

\begin{figure}[h]\centering
\begin{minipage}[t]{0.49\textwidth}
\centering
\includegraphics[width=\textwidth]{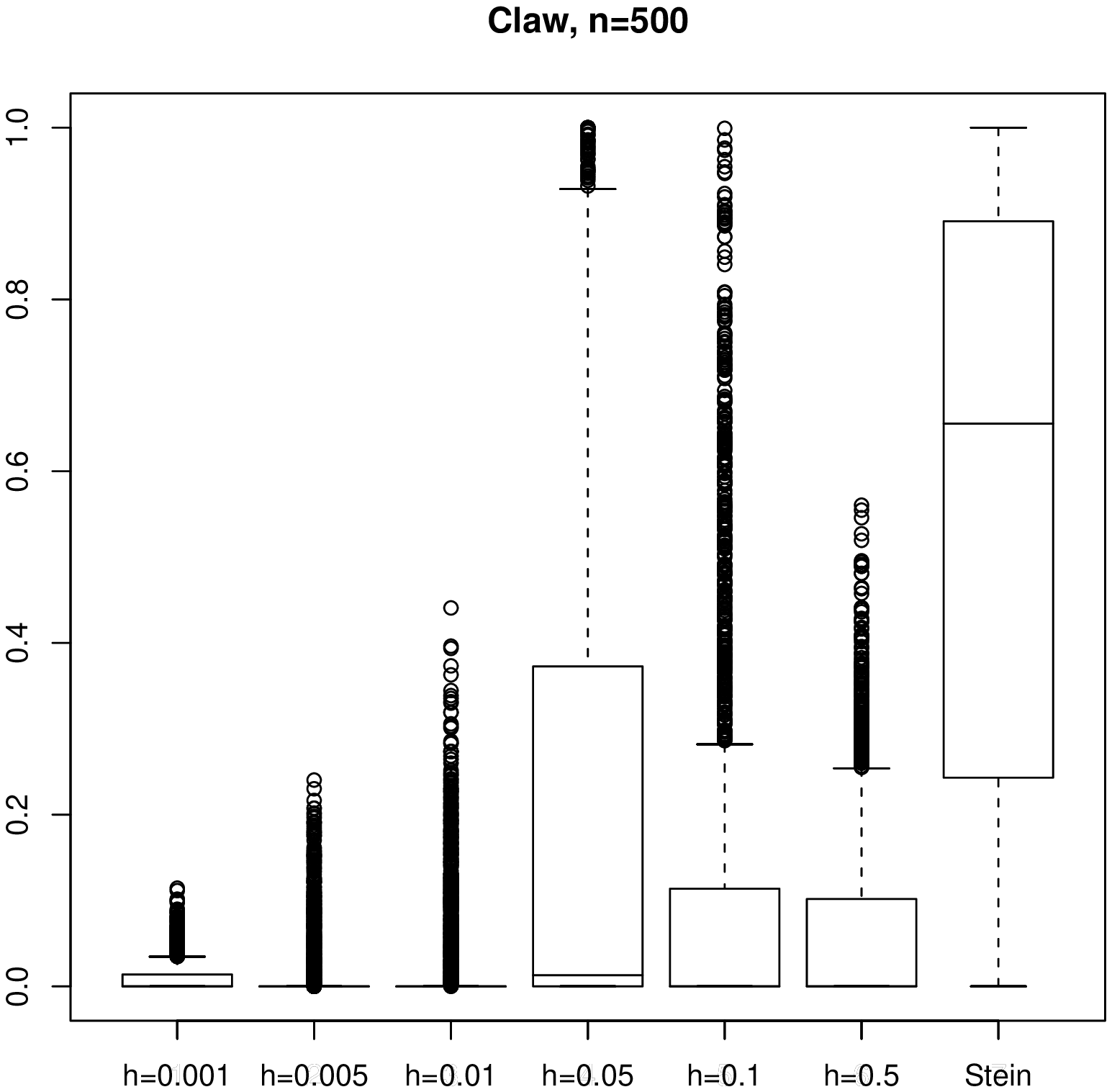}
\end{minipage}
\begin{minipage}[t]{0.49\textwidth}
\centering
\includegraphics[width=\textwidth]{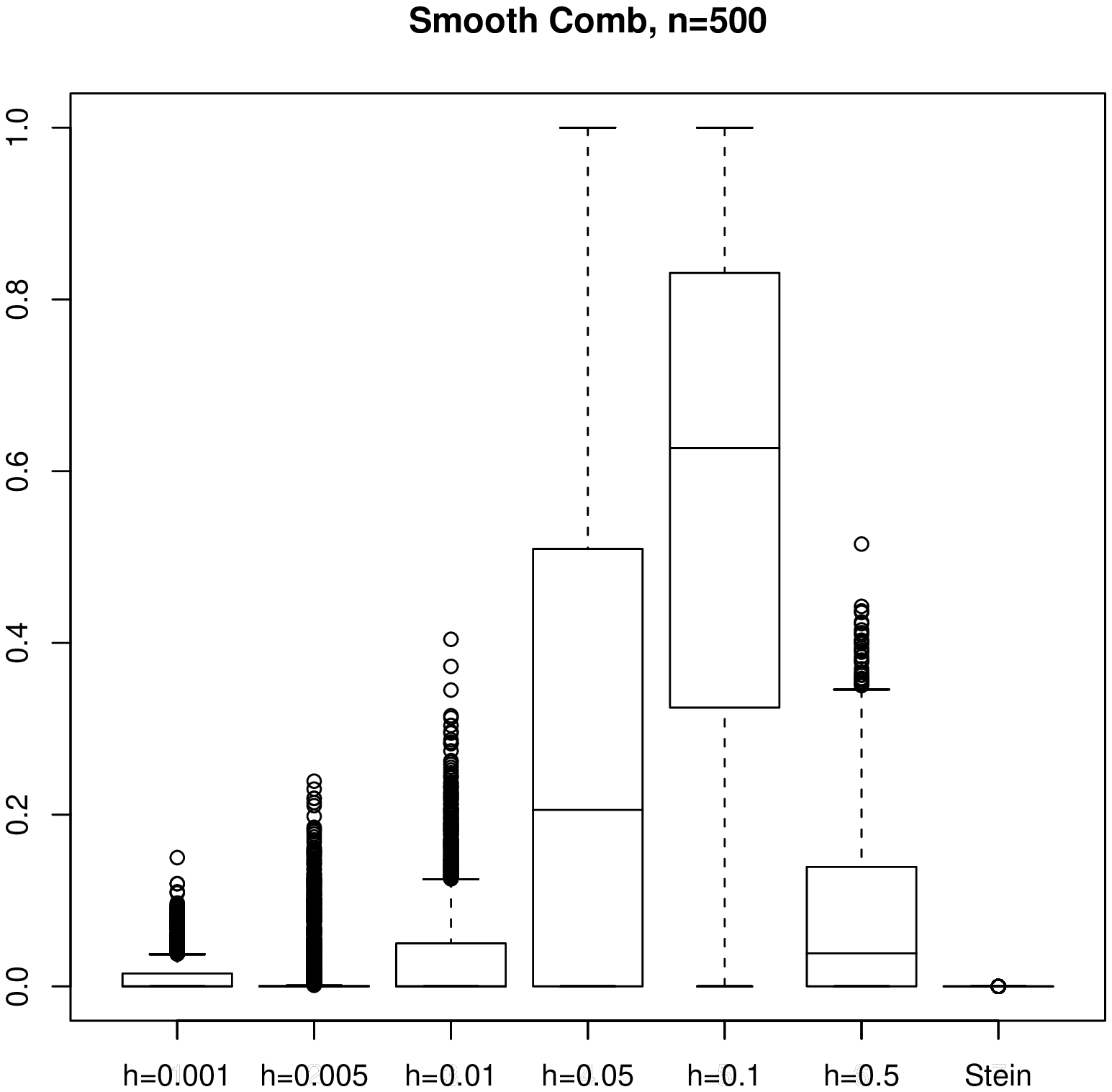}
\end{minipage}
\caption{Boxplots for the Claw and Smooth Comb densities
}\label{boxp}
\end{figure}

{\tt AggStein} is slightly better than {\tt AggPure} for the Claw
density and outperforms the other tested estimators in almost all
the considered cases, so we studied this procedure in more detail.
We focused on the Claw and Smooth Comb densities and a sample of
size 500. Figure~\ref{compag} gives a visual comparison of the {\tt
AggStein} procedure and the DPI procedure.
It illustrates the oversmoothing effect of the DPI procedure and the
fact that the {\tt AggStein} procedure adapts to inhomogeneous
smoothness. We finally comment on two other aspects of the {\tt
AggStein} procedure:
\begin{itemize}
\item the distribution of weights that are allocated to the aggregated
estimators,
\item the robustness to the number and size of the splits.
\end{itemize}
The boxplots represented in Figure~\ref{boxp} give the distributions
of weights allocated to 7 estimators to be aggregated, the 6 kernel
density estimators and the block Stein estimator. The boxplots are
constructed from 2000 values of the vector of the weights (200
samples times 10 splits).
We immediately notice that for the Claw density a median weight
greater than 0.65 is allocated to the block Stein estimator. This
can be explained by the fact that the block Stein estimator performs
better than kernel density estimators on this density [cf. MISE of
{\tt AggPure} and Stein in Table~\ref{misecsc} (left)], and the {\tt
AggStein} procedure takes advantage of it. On the other hand, for
the Smooth Comb density, the block Stein estimator does not perform
significantly better than the kernel density estimators [see
Table~\ref{misecsc} (right)] and the {\tt AggStein} procedure does
not use it at all. For this sample size and this density, the
procedures {\tt AggStein} and {\tt AggPure} are equivalent.

A free parameter of the aggregation procedures is the set of splits.
In this study we choose random splits and we only have to specify
their number and sizes. Obviously, we are interested to have less
splits in order to make the procedure less time consuming. Figure
\ref{sensplit} gives the sensibility of MISE both to the number of
splits and to the size of the training sample in the case of dens1
and dens2 with the overall sample size 200.
Two important conclusions are: (i) there exists a size of the
training sample that achieves the minimum MISE, and (ii) there is
essentially nothing to gain by producing more than 20 splits.
Similar results are obtained for {\tt AggPure}, and they are valid
on the whole set of tested densities.


{\bf Acknowledgment:} We would like to thank the referees for
helpful remarks and Lucien Birg\'e for suggesting an improvement of
the constants in Theorem~\ref{t1l} as well as a simplification of
its proof. We refer to Birg\'e (2006) for comments on a previous
version of this paper.

\nopagebreak

\end{document}